\documentclass[letterpaper]{article}

\usepackage{arxiv}

\usepackage{style}

\usepackage{lineno,hyperref}

\usepackage{amsmath, amssymb, latexsym, multicol, amsthm, bm, placeins}
\usepackage[export]{adjustbox}
\usepackage{color}
\usepackage{multirow}

\usepackage{style}

\makeatletter\@addtoreset{equation}{section}\makeatother

\renewcommand{\shorttitle}{Higher-order adaptive virtual element methods with contraction properties}

\numberwithin{equation}{section}

\makeatletter
\renewcommand{\@biblabel}[1]{#1\hfill \hspace{-0.2cm}}
\makeatother

\begin{document}

\title{Higher-order adaptive virtual element methods with contraction properties}

\author{%
  Claudio Canuto\affil{1,}\corrauth
  and
  Davide Fassino\affil{1}
}

\author{C. Canuto\thanks{Department of Mathematical Sciences “Giuseppe Luigi Lagrange”, Politecnico di Torino, Corso Duca degli Abruzzi 24, Torino, 10129, Italy. claudio.canuto@polito.it (C. Canuto), davide.fassino@polito.it (D. Fassino).}
	\And
	D. Fassino\footnotemark[1]
}

\maketitle
\begin{abstract} 

The realization of a standard Adaptive Finite Element Method (AFEM) preserves the  mesh conformity by performing a completion step in the refinement loop:
in addition to elements marked for refinement due to their contribution
to the global error estimator, other elements are refined.

In the new perspective opened by the introduction of Virtual Element
Methods (VEM), elements with hanging nodes can be
viewed as polygons with aligned edges, carrying virtual functions together with standard polynomial functions.
The potential advantage is that all activated degrees of freedom are
motivated by error reduction, not just by geometric reasons. 

This point of view is at the basis of the paper [L. Beir\~{a}o da Veiga et al., “Adaptive VEM: stabilization-free a posteriori error analysis and contraction
property”, SIAM Journal on Numerical Analysis, vol. 61,
2023], devoted to the convergence analysis of an adaptive VEM
generated by the successive newest-vertex bisections of triangular
elements without applying completion, in the lowest-order case
(polynomial degree $k=1$).

The purpose of this paper is to extend these results to the case of VEMs of order $k\ge2$ built on
triangular meshes. The problem at hand is a variable-coefficient,
second-order self-adjoint elliptic equation with Dirichlet boundary
conditions; the data of the problem are assumed to be piecewise polynomials of degree $k-1$. By extending the concept of global index of a hanging node, under an admissibility assumption of the mesh, we derive a stabilization-free a posteriori error estimator. This is the sum of residual-type terms and certain virtual inconsistency terms (which vanish for $k=1$). We define an adaptive VEM of order $k$ based on this estimator, and we prove its convergence by establishing a contraction result for a linear combination of
(squared) energy norm of the error,  (squared) residual estimator, and
(squared) virtual inconsistency estimator.

\end{abstract}


\maketitle
\section{Introduction}
Adaptive Finite Element Methods (AFEM) for self-adjoint coercive
problems written in the form
\begin{linenomath}
\begin{equation*}
u \in V \ : \ \mathcal{B}(u,v) = F(v) \qquad \forall v \in \V
\end{equation*}\end{linenomath}
iterate the sequence $ {\tt SOLVE} \rightarrow {\tt ESTIMATE}
\rightarrow {\tt MARK} \rightarrow {\tt REFINE}$ to produce better and
better approximations of $u$. Their practical efficiency is corroborated
by sound theoretical results of convergence, complexity, and optimality,
which in various cases (such as, e.g., conforming $h$-versions)
completely explain the behaviour of the adaptive algorithms \cite{Dorfler, BDD:04, CasconNochetto, VeeserNochetto, AxiomsOfAdaptivity}.

The standard AFEM realization preserves the conformity of the initial
mesh, at the expense of performing a completion step in ${\tt REFINE}$:
in addition to elements marked for refinement due to their contribution
to the global error estimator, other elements are refined. Without this
step, one would obtain nonconforming meshes, containing elements with
hanging nodes.

In the new perspective opened by the introduction of Virtual Element
Methods (VEM) \cite{BasicPrinciples, TheHitchhikersGuide}, elements with hanging nodes can be
viewed as polygons with aligned edges, carrying virtual (i.e.,
non-accessible) functions together with standard polynomial functions.
The potential advantage is that all activated degrees of freedom are
motivated by error reduction, not just by geometric reasons. On the
other hand, in this transformation of an adaptive FEM into an adaptive
VEM, one looses the availability of a general convergence theory, which
so far is lacking (although results on a posteriori error estimates
\cite{BeiraoManzini,Cangiani} have been obtained, together with efficient practical recipes
for refining polytopal meshes \cite{BerroneBorio, BerroneDauria, AntoniettiDassiManuzzi}).

Such a shift in perspective inspired the recent papers
\cite{AVEMstabfree,AVEMConvergenceOptimality}, devoted to the analysis of an adaptive VEM
generated by the successive newest-vertex bisections of triangular
elements without applying completion, in the lowest-order case
(polynomial degree $k=1$). Despite the simple geometric setup, the
investigation faced  some VEM-specific obstacles in the analysis, giving
answers that could prove useful in the study of more general adaptive
VEM discretizations. For instance, a VEM solution $\umesh \in \Vmesh
\subset \V$, defined by the Galerkin projection
\begin{linenomath}
\begin{equation*}
u_\mesh \in \Vmesh \ : \ \mathcal{B}_\mesh(\umesh,\vmesh) =
F_\mesh(\vmesh) \qquad \forall \vmesh \in \Vmesh \,,
\end{equation*}
\end{linenomath}
satisfies an a posteriori error bound of the type
\begin{linenomath}
\begin{equation*}
\Vert u - \umesh \Vert_\V^2 \ \lesssim\  \etamesh^2(\umesh) +
S_\mesh(\umesh,\umesh) \,,
\end{equation*}
\end{linenomath}
where $\etamesh(\umesh)$ is a residual-type error estimator,
$S_\mesh(\umesh,\umesh)$ is the stabilization term that makes the
discrete bilinear form $ \mathcal{B}_\mesh(\umesh,\vmesh)$ coercive in
$\V$, and for simplicity we assume piecewise constant data on the mesh
$\mesh$. Unfortunately, the term $S_\mesh(\umesh,\umesh)$ need not
reduce under a mesh refinement, as $\etamesh^2(\umesh)$ does: this makes
the convergence analysis problematic. However, one of the key results
obtained in \cite{AVEMstabfree} states that
$S_\mesh(\umesh,\umesh)$ is dominated by $\etamesh^2(\umesh)$, i.e.
\begin{linenomath}
\begin{equation*}
S_\mesh(\umesh,\umesh) \ \lesssim\  \etamesh^2(\umesh)  \,,
\end{equation*}\end{linenomath}
provided an assumption of {\em admissibility} of the non-conforming
meshes generated by successive refinements is fulfilled; such a
restriction, which appears to have little practical impact, amounts to
requiring the uniform boundedness of the {\em global index} of all
hanging node, a useful concept introduced in \cite{AVEMstabfree} to
hierarchically  organize the set of hanging nodes. Once the a posteriori
error bound is reduced to
\begin{linenomath}
\begin{equation*}
\Vert u - \umesh \Vert_\V^2 \ \lesssim\  \etamesh^2(\umesh)  \,,
\end{equation*}\end{linenomath}
the convergence analysis becomes feasible, and a contraction property is
proven to hold for a  linear combination of the (squared) energy norm of
the error and the (squared) residual estimator.

The purpose of this paper is to extend the results in
\cite{AVEMstabfree} to the case of VEMs of order $k\ge2$ built on
triangular meshes. Note that the interest in avoiding the creation of new elements just to satisfy the conformity condition of the mesh  becomes more and more evident as the polynomial degree increases. 
The problem at hand is again a variable-coefficient,
second-order self-adjoint elliptic equation with Dirichlet boundary
conditions. The geometric concept of hanging node (a vertex for some
elements, contained inside an edge of some other elements) is replaced
by a functional one, referring to the degrees of freedom associated with
the node; once the meaning of hanging node is clarified, the definition
of {\em global index} of a node, and its role in the analysis, is
similar to the one given in \cite{AVEMstabfree}.

A significant difference with respect to the content of that paper
concerns the control of the stabilization term, which does not involve
only the residual estimator, but a new term, called the {\em virtual
inconsistency estimator} and denoted by $\Psi_\mesh(\umesh)$. It
measures the projection error, upon local spaces of polynomials, of
certain expressions depending on the operator coefficients and the
discrete solution; it vanishes when $k=1$ or when the coefficients are
constant. The new stabilization bound, which we derive under an
admissibility assumption of the mesh,  takes the form
\begin{linenomath}
\begin{equation*}
S_\mesh(\umesh,\umesh) \ \lesssim\  \etamesh^2(\umesh)  +
\Psi_\mesh^2(\umesh)  \,,
\end{equation*}\end{linenomath}
which leads to the a posteriori, stabilization-free error control
\begin{linenomath}
\begin{equation*}
\Vert u - \umesh \Vert_\V^2 \ \lesssim\  \etamesh^2(\umesh) +
\Psi_\mesh^2(\umesh)   \,.
\end{equation*}\end{linenomath}
Correspondingly, we obtain the convergence of the adaptive VEM of order
$k$ by proving a contraction result for a linear combination of
(squared) energy norm of the error,  (squared) residual estimator, and
(squared) virtual inconsistency estimator.

Similarly to \cite{AVEMstabfree}, we assume here that the data $\data$
of our boundary-value problem are piecewise polynomials of degrees
related to $k-1$, on the initial mesh $\mesh_0$ and consequently on each
mesh $\mesh$ derived by newest-vertex bisection. This is not a
restriction, since we propose to insert the adaptive VEM procedure just
described, which we now consider as a module ${\tt GALERKIN}$, into an
outer loop ${\tt AVEM}$ of the form
\begin{algorithmic}
\State$[\mathcal{T}, u_\mathcal{T}] ={\tt
AVEM}(\mathcal{T}_0,\epsilon_0,\omega, \text{tol})$
     \State $j=0$
     \While{$\epsilon_j > \frac{1}{2} \mathrm{tol}$}
      \State   $[\hat{\mathcal{T}}_j, \hat{\mathcal{D}}_j] = {\tt
DATA}(\mathcal{T}_j, \mathcal{D},  \epsilon_j)$
        \State $[{\mathcal{T}}_{j+1}, \data_{j+1}] = {\tt
GALERKIN}(\hat{\mathcal{T}}_j, \hat{\mathcal{D}}_j, \epsilon_j)$
        \State  $\epsilon_{j+1} \gets\frac{1}{2}\epsilon_j$
        \State $j \gets j +1$
     \EndWhile\\
     \Return
\end{algorithmic}
where the module ${\tt DATA}$ produces, via greedy-type iterations, a
piecewise polynomial approximation of the input data with prescribed
accuracy, defined on a suitable refinement of the input partition.
Manifestly, the target accuracy is matched after a finite number of
calls to ${\tt DATA}$ and ${\tt GALERKIN}$.
Properties of complexity and quasi-optimality of this two-loop algorithm
are investigated in \cite{AVEMConvergenceOptimality} in the linear case $k=1$. We plan to do
the same for the case $k \geq 2$ in a forthcoming paper.

The outline of this paper is as follows. In Sections \ref{sec:VEMspaces}
and \ref{sec:discretization}, we introduce the model boundary-value
problem, and its discretization by an enhanced version of the VEM
(\cite{EquivalentProjectors}).  In Section \ref{sec:indexofanode} we define the global
index of a node, and we formulate the admissibility assumption on the
mesh. Two essential properties for bounding the stabilization term are
established  in Section \ref{sec:TheKeyProperties}.  The a posteriori
error estimators are defined in Section
\ref{sec:aPosterioriErrorEstimator},  whereas stabilization-free a
posteriori error estimates are proven in Section
\ref{sec:aPosterioriErrorEstimates}. In Section
\ref{sec:EffectOfAMeshRefinement}, we investigate how the a posteriori
error estimators are reduced under mesh refinement. These properties are
needed to justify the refinement strategy in our adaptive module {\tt
GALERKIN}, which is described in Section \ref{Sec:GALERKIN}. The paper
ends with the proof of convergence of the loop {\tt GALERKIN}, reported
in Section \ref{sec:galerkin}.


\section{VEM spaces of order $k\ge2$}\label{sec:VEMspaces}

We consider the following Dirichlet boundary value problem in a  polygonal domain $\Omega$,
\begin{linenomath}
\begin{align}
    \begin{cases}
         -\nabla \cdot (A \nabla u) + c u = f &\text{ in $\Omega$,}\\
         u = 0 &\text{ on $\partial \Omega$,}
    \end{cases}
    \label{Initial_Problem}
\end{align}
\end{linenomath}
where $A \in (L^\infty(\Omega))^{2\times 2}$ is symmetric and uniformly positive definite in $\Omega$, $c \in L^{\infty}(\Omega)$ and non-negative in $\Omega$, $f \in L^2 (\Omega)$. Data will be denoted by $\data = (A,c,f)$. The variational formulation of this problem is written as
\begin{linenomath}
\begin{align}
\begin{cases}
\text{find }u \in \mathbb{V}  :=H^1_{0}(\Omega) &\text{ such that}\\
\mathcal{B}(u,v) = (f,v), &\forall \;v \in \mathbb{V},
\end{cases}
\label{Variational_Problem}
\end{align}
\end{linenomath}
where $(\cdot,\cdot)$ is the scalar product in $L^2(\Omega)$ and $\mathcal{B}(u,v) := a(u,v) + m(u,v) $ is the bilinear form associated with Problem \eqref{Initial_Problem}, i.e,
\begin{linenomath}
  \begin{align*}
    &a(u,v) := (A\; \nabla u , \nabla v) & m(u,v) := ( c\; u,v).
\end{align*}  
\end{linenomath}
We denote the energy norm as $\EnergyNorm{\cdot}= \sqrt{\mathcal{B}(\cdot,\cdot)}$, which satisfies
\begin{linenomath}
\begin{align}\label{bound_energynorm}
    &c_\mathcal{B} \normH{v}^2\le \EnergyNorm{v}^2 \le c^\mathcal{B}\normH{v}^2, &\forall v \in \mathbb{V},
\end{align}
\end{linenomath}
for suitable $0 < c_\mathcal{B} \le c^\mathcal{B}$.

In order to find a discrete approximation of the solution of Problem \eqref{Variational_Problem}, we firstly introduce a fixed initial partition $\mathcal{T}_0$ on the domain $\overline{\Omega}$ made of triangular elements $E$. We will denote by $\mathcal{T}$ any refinement of $\mathcal{T}_0$ obtained by a finite number of newest-vertex element bisections. We underline that we are not requiring $\mathcal{T}$ to be a conforming mesh, since hanging nodes may arise in the refinement. The classification of nodes, which will play a crucial role in the proofs presented in this paper, is postponed in Section \ref{sec:indexofanode}.

According to the Virtual Element theory \cite{BasicPrinciples}, an element $E$ of the triangulation can be viewed as a polygon with more than three edges, if some hanging nodes are sitting on its boundary. We can then denote by $\mathcal{E}_E$ the set of edges $e$ of element $E$ and $\mathcal{E}:= \bigcup_{E\in\mathcal{T}}\mathcal{E}_E$. We finally define the diameter of an element $E$ as $h_E= |E|^{1/2}$ and $h=\max_{E\in \mathcal{T}}\{h_E\}$.

We introduce the functional spaces needed to apply the Virtual Element Method (VEM). We start by defining the space of functions on the boundary of $E$, $\mathbb{V}_{\partial E,k}$, which is constituted by the functions that are continuous on the boundary of $E$ and that, when restricted to any edge of $\partial E$, are polynomials of degree $k>0$, i.e, 
\begin{linenomath}
\begin{align*}
    \mathbb{V}_{\partial E,k}:= \{v\in C^{0}\left (\partial E\right): v\vert_e \in \mathbb{P}_k(e), \forall e \subset \partial E\}.
\end{align*}
\end{linenomath}
Then, we define the ``enhanced'' VEM space in $E$, as done in \cite{EquivalentProjectors}, such that
\begin{linenomath}
\begin{align}\label{Eq:enhanced}
     \mathbb{V}_{E,k}:= \left\{v\in H^1\left (E\right)\;:\; v\vert_{\partial E} \in  \mathbb{V}_{\partial E,k},\; \Delta v \in  \mathbb{P}_{k}(E),\;(v - \Pi^\nabla_E v, q)_E=0\; \forall q \in \mathbb{P}_k(E) \right\}, 
\end{align}
\end{linenomath}
where $\Pi^\nabla_{E} \;:\;H^1(E)\rightarrow \mathbb{P}_k(E)$ is the projector defined by
\begin{linenomath}
\begin{align*}
    &(\nabla (v - \Pi^\nabla_E v), \nabla q)_E = 0\;\;\; \forall q \in \mathbb{P}_k(E), &\int_{\partial E}(v -\Pi^\nabla_E v)=0. 
\end{align*}
\end{linenomath}
We remark that $\mathbb{V}_{E,k}$ contains the polynomial space of degree $k$ on $E$ and its dimension is
\begin{linenomath}
\begin{align}\label{dimVE}
    \dim(\mathbb{V}_{E,k})= 3 k + \frac{k (k-1)}{2},
\end{align}
\end{linenomath}
since in our analysis we consider triangular elements. We notice that in the case $k>1$ a function $v$ in $\mathbb{V}_{E,k}$ is uniquely defined by 
\begin{itemize}
    \item the set of the values at the vertices of $E$;
    \item the set of the values at the $k-1$ equally-spaced internal points on each edge of $\partial E$;
    \item  the set of the moments
    $\frac{1}{|E|}\int_{E} v(\bm{x}) m(\bm{x}) d\bm{x}$ $\forall m \in \mathcal{M}_{k-2}(E)$,
    \end{itemize}
where the set $\mathcal{M}_{p}(E)$, $p \geq 0$, is defined as
\begin{linenomath}
\begin{align}\label{ME}
    \mathcal{M}_{p}(E) = \left\{\left( \frac{\bm{x}-\bm{x}_E}{h_E}\right)^s, |s|\le p\right\}.
\end{align}
\end{linenomath}
We will denote by $\bm{\mu}_p(E,v) = \left(\frac{1}{|E|}\int_{E} v(\bm{x}) m(\bm{x}) d\bm{x} : m \in \mathcal{M}_{p}(E) \setminus \mathcal{M}_{p-1}(E) \right)$ the vector of the moments of $v$ of order $p$. By $|\bm{\mu}_p(E,v) |$ we will denote the $l^2$-norm of this vector.

\medskip
We can now introduce the global discrete space as 
\begin{linenomath}
  \begin{align*}
    \mathbb{V}_{\mathcal{T}}:= \{ v \in \mathbb{V}: \; v|_E \in \mathbb{V}_{E, k}\; \forall E \in \mathcal{T} \}.
\end{align*}  
\end{linenomath}
On $\mathcal{T}$ we need also to give the definition of the space of piecewise polynomial functions on $\mathcal{T}$
\begin{linenomath}
\begin{align}\label{def:WkT}
    \mathbb{W}^k_{\mathcal{T}}:= \{ w \in L^2(\Omega): w|_E \in \mathbb{P}_k(E)\;\forall E \in \mathcal{T}\},
\end{align}
\end{linenomath}
and its subspace
\begin{linenomath}
\begin{align}\label{def:V0T}
    \mathbb{V}^0_{\mathcal{T}}:= \mathbb{V}_{\mathcal{T}} \cap \mathbb{W}^k_{\mathcal{T}},
\end{align}
\end{linenomath}
which plays a crucial role in the forthcoming analysis.

We now introduce a series of projectors that will be used in the rest of the paper. For any $E\in \mathcal{T}$, we denote by $\Pi^0_{p,E}: L^2(E)\rightarrow \mathbb{P}_p(E)$ the $L^2(E)$-orthogonal projector onto the space of polynomial of degree $p$ on $E$.  Thanks to the choice of the enhanced space $\mathbb{V}_{E,k}$ \eqref{Eq:enhanced}, we remark that $\Pi^0_{k,E} v$ and $\Pi^0_{k-1, E}\nabla v$ can be computed for any function $v \in \mathbb{V}_{E,k}$, see \cite{EquivalentProjectors} for the details. To simplify the notation, in the following we will drop the symbol $E$ from $\Pi^0_{k, E}$ when no confusion arises. The global $L^2$-orthogonal projector is denoted by $\Pi^0_{p,\mathcal{T}}: L^2(\Omega)\rightarrow \mathbb{W}^p_\mathcal{T}$.

We can also define the Lagrange interpolation operator $\mathcal{I}_E:\mathbb{V}_{E,k}\rightarrow\mathbb{P}_k(E)$ on $E$, which builds a polynomial of  degree $k$ using the $3 k$ degrees of freedom on the boundary of $E$ and the moments of order $\leq k-3$, since
\begin{linenomath}
\begin{align*}
 \dim(\mathbb{P}_k(E)) = 3k + \frac{(k-1)(k-2)}{2}.
\end{align*}
\end{linenomath}
Moreover, we will denote by $\mathcal{I}_\mathcal{T}:\mathbb{V}_\mathcal{T}\rightarrow \mathbb{W}^k_{\mathcal{T}}$ the Lagrange interpolation operator that restricts to $\mathcal{I}_E$ on each $E\in\mathcal{T}$.

\section{Discretization with data of degree $k-1$}\label{sec:discretization}

In the rest of this paper, we assume that data $\data=(A,c,f)$ are piecewise polynomials of degree $k-1$ on the initial partition $\mesh_0$, hence on each partition $\mesh$ obtained by newest-vertex refinement. Their values on each element of the triangulation will be denoted by $(A_E,c_E, f_E) \in (\mathbb{P}_{k-1}(E))^{2 \times 2} \times \mathbb{P}_{k-1}(E) \times \mathbb{P}_{k-1}(E)$.

We here define the bilinear forms that we need for the Galerkin discretization problem, starting from $a_E, m_E : \mathbb{V}_{E,k} \times \mathbb{V}_{E,k} \rightarrow \mathbb{R}$, such that
\begin{linenomath}
    \begin{align*}
        &a_\mathcal{T}(v,w) := \sum_{E \in \mathcal{T}}\int_E A_E \left(\Pi^0_{k-1} \nabla v\right) \left(\Pi^0_{k-1} \nabla w\right)=:\sum_{E\in \mathcal{T}}a_E(v,w),\\ & m_\mathcal{T}(v,w) := \sum_{E \in \mathcal{T}} \int_E c_E\;\Pi ^0_{k} v \;\Pi^0_{k} w =:\sum_{E\in \mathcal{T}}m_E (v, w).
    \end{align*}
\end{linenomath}
We also introduce the symmetric bilinear form $s_E: \mathbb{V}_E \times \mathbb{V}_E \rightarrow \mathbb{R}$ as
\begin{linenomath}
\begin{align*}
 s_E (v,w) := \sum^{\overline{\mathcal{N}}_E}_{i=1} v(\bm{x}_i) w(\bm{x}_i),
\end{align*}
\end{linenomath}
where $\{\bm{x}_i\}^{\overline{\mathcal{N}}_E}_{i=1}$ indicates the set of the degrees of freedom on the boundary of $E$. Indeed, we remark that in this case the stabilization term can be built without using the internal degrees of freedom, as shown in \cite{StabTerm}. We assume for $s_E$ the existence of two positive constant $c_s$ and $C_s$ independent on $E$, such that
\begin{linenomath}
\begin{align}\label{StabilityStabilizationForm}
    &c_s\normHE{v}^2 \le s_E(v,v) \le C_s \normHE{v}^2 &\forall v, w \in \mathbb{V}_E\setminus \mathbb{R}.
\end{align}
\end{linenomath}
We define the local stabilizing form as
\begin{linenomath}
    \begin{align*}
    &S_E(v, w)= s_E(v - \mathcal{I}_E v, w - \mathcal{I}_E w)
    &\forall v, w \in \mathbb{V}_E,
\end{align*}
\end{linenomath}
and the global stabilization form
\begin{linenomath}
\begin{align*}
 &S_\mathcal{T} (v,w):= \sum_{E\in \mathcal{T}}S_E(v,w) &\forall v, w \in \mathbb{V}_\mathcal{T}.
\end{align*}
\end{linenomath}
From \eqref{StabilityStabilizationForm}, we get
\begin{linenomath}
\begin{align*}
    &S_\mathcal{T}(v,v) \simeq \normHT{v- \mathcal{I}_\mathcal{T}v}^2 &\forall v \in \mathbb{V}_\mathcal{T}, 
\end{align*}
\end{linenomath}
where $\normHT{\, \cdot \, }$ denotes the broken $H^1$-seminorm over $\mathcal{T}$.
Thus, we can now define the bilinear form $\mathcal{B}_{\mathcal{T}}(\cdot,\cdot)$, $\mathcal{B}_{\mathcal{T}}: \mathbb{V}_{\mathcal{T}} \times \mathbb{V}_{\mathcal{T}}\rightarrow \mathbb{R}$, as
\begin{linenomath}
    \begin{align}\label{DefBT}
        \mathcal{B}_\mathcal{T}(v,w)= a_\mathcal{T}(v,w)+ m_\mathcal{T}(v,w) +\gamma S_\mathcal{T}(v, w),
    \end{align}
\end{linenomath}
with $\gamma$ independent of $\mesh$ satisfying $\gamma \ge \gamma_0$ for some fixed $\gamma_0>0$. 
For the loading term we introduce $\mathcal{F}_\mathcal{T}: \mathbb{V}_\mathcal{T} \rightarrow \mathbb{R}$ as
\begin{linenomath}
\begin{align}\label{Def:forcing_term}
 &\mathcal{F}_\mathcal{T}(v):= \sum_{E\in \mathcal{T}} \int_E f_E\; \Pi^0_k v  = \sum_{E\in \mathcal{T}} \int_E f_E v &\forall v \in \mathbb{V}_\mathcal{T} \,,
\end{align}
\end{linenomath}
since $f_E$ has been already approximated with a polynomial of degree $k-1$. Note that the equality in \eqref{Def:forcing_term} remains true if $f_E$ is an approximation of $f$ of degree $k$ on $E$.

We have now defined all the forms that appear in the discrete formulation of the Problem \eqref{Variational_Problem}. It reads as
\begin{linenomath}
\begin{align}
\begin{cases}
\text{find }u_{\mathcal{T}} \; \in \; \mathbb{V}_{\mathcal{T}} \text{ such that}\\
\mathcal{B}_\mathcal{T}(u_{\mathcal{T}},v) = \mathcal{F}_\mathcal{T}(v),\; \forall \;v \in \mathbb{V}_{\mathcal{T}}.
\end{cases}
\label{Discrete_Variazional_Problem}
\end{align}
\end{linenomath}
The bilinear form $\mathcal{B}_\mathcal{T}$ is continuous and coercive, hence, there exists a unique and stable solution of the Problem \eqref{Discrete_Variazional_Problem}.
Furthermore, the following result extends Lemma 2.6 in \cite{AVEMstabfree}.
\begin{lemma}[Gakerkin quasi-orthogonality]\label{Lemma:QuasiOrtog}
    For any $v\in \mathbb{V}_\mathcal{T}$ and $w \in \mathbb{V}^0_\mathcal{T}$, it holds
    \begin{linenomath}    
    \begin{equation*}
    \begin{split}
        a_\mathcal{T}(v, w)&=a(v, w) - \sum_{E \in \mesh} \int_E A_E \left((I - \Pi^0_{k-1}) \nabla v\right)  \nabla w, \\
        m_\mathcal{T}(v, w)&=m(v, w) - \sum_{E \in \mesh} \int_E c_E \left((I-\Pi^0_{k}) v\right)  w, \\
        S_\mathcal{T}(v,w)&=0\,.
  \end{split}
  \end{equation*}
  \end{linenomath}
  Consequently,
  \begin{linenomath}
    \begin{align*}
        | \mathcal{B}(u - u_\mathcal{T}, w)| \lesssim S_\mathcal{T}(u_\mesh, u_\mesh)^{1/2} |w|_{1,\Omega} \,,
    \end{align*}
    \end{linenomath}
    where $u$ is the solution of \eqref{Variational_Problem} and $u_\mathcal{T}$ the solution of \eqref{Discrete_Variazional_Problem}.
\end{lemma}

\section{The index of a node}\label{sec:indexofanode}

A crucial concept, firstly introduced in \cite{AVEMstabfree} for the case $k=1$, is the \textit{global index of a node}: it will be used in the proofs of Section \ref{sec:TheKeyProperties}. In order to extend its definition to the case $k >1$, we preliminarily introduce some useful definitions.

Let $\hat{E}:= \{(x,y)\in \mathbb{R}^2:\; x\ge 0,\; y\ge 0,\; x+y\le 1\}$ be the reference element
and denote by $\hat{R}_{\hat{E},k}$ the $k$-lattice built on $\hat{E}$, i.e.,
\begin{linenomath}
\begin{align*}
    \hat{R}_{\hat{E},k}:= \left\{\left(\frac{i}{k},\frac{j}{k}\right)\in \mathbb{R}^2:\; i\ge 0,\; j\ge 0,\; i+j\le k\right\}.
\end{align*}
\end{linenomath}
Considering the affine function $F_E: \hat{E}\rightarrow E$ mapping the reference element onto an element $E\in \mesh$, we define the physical lattice on $E$ by
    $R_{E, k}:= F_E(\hat{R}_{\hat{E},k})$,
and the set of \textit{proper nodes} of $E$ as the points of the physical lattice sitting on the boundary of $E$, i.e.,
\begin{linenomath}
\begin{align*}
  \mathcal{P}_E:= R_{E, k}\cap \partial E.  
\end{align*}
\end{linenomath}
Observe that we implicitly assume that $k\ge2$ is sufficiently small so that interpolation on equally spaced nodes is numerically stable.

Next, we denote by $\mathcal{H}_E$ the set of \textit{hanging nodes} of $E$, i.e., the set of points $\bm{x} \in \partial E$ that are not proper nodes of $E$, but that are proper nodes of some other element $E'$, i.e.,
\begin{linenomath}
\begin{align*}
    \mathcal{H}_E:= \left \{ \bm{x} \in \partial E: \exists E'  \in \mesh \text{ such that } \bm{x}\in \mathcal{P}_{E'}\right  \} \setminus \mathcal{P}_E.
\end{align*}
\end{linenomath}
Finally, let $\mathcal{N}_E:= \mathcal{P}_E \cup \mathcal{H}_E$ be the set of all nodes sitting on  $E$. 

At the global level, $\mathcal{N} :=  \bigcup_{E \in\mathcal{T}} \mathcal{N}_E$  will be the set of all nodes of the triangulation $\mathcal{T}$, which we split into the set $\mathcal{P}:= \left\{\bm{x}\in \mathcal{N}: \bm{x}\in \mathcal{P}_E \ \forall E \text{ containing }\bm{x}\right\}$ of the \textit{proper nodes of }$\mesh$,  and the set $\mathcal{H}:= \mathcal{N}\setminus \mathcal{P}$ of the \textit{hanging nodes of }$\mesh$.

\begin{figure}[t!]
    \centering
    \subfloat[]{
\begin{tikzpicture}
\draw [thick,blue] plot coordinates{(8,0)  (0,0)};
\draw[very thick, blue] plot[mark =square, only marks, mark size=3pt] coordinates{(0,0)(2,0)(4,0)(6,0)(8,0)};
\coordinate [label= {[blue]$S$}] (z1) at (4.1,0.3);
\coordinate [label= {$\bm{\xi}_1$}] (z1) at (0,-0.7);
\coordinate [label= {$\bm{\xi}_{k+1}$}] (z1) at (8,-0.7);
\end{tikzpicture}
}\\
\centering
\subfloat[]{
\begin{tikzpicture}
\draw [thick] plot coordinates{(4,0)  (0,0)};
\draw [thick,orange] plot coordinates{(4,0)  (8,0)};
\draw[very thick,blue] plot[mark =square, only marks, mark size=3pt] coordinates{(0,0)(2,0)(4,0)(6,0)(8,0)};
\draw[very thick, red] plot[mark =*, only marks, mark size=2pt] coordinates{(0,0)(1,0)(2,0)(3,0)(4,0)(5,0)(6,0)(7,0)(8,0)};
\coordinate [label=  { [blue]$S$}](z1) at (4.1,0.3);
\coordinate [label= $S^-$] (z1) at (1.9,-0.7);
\coordinate [label= { [orange]$S^+$}] (z1) at (6.1,-0.7);
\coordinate [label= {$\bm{\zeta}_1$}] (z1) at (1,0.1);
\coordinate [label= {$\bm{\zeta}_{k}$}] (z1) at (7,0.1);
\end{tikzpicture}}
\caption{Blue squares represent the $k+1$ equally-spaced nodes $\bm{\xi}_n$ on the edge $S$ before refinement. Red circles represent the $2k +1$ nodes that arise after  refinement. We have denoted by $\bm{\zeta}_i$  the new nodes that do not coincide with any $\bm{\xi}_n$.}
\label{Figure:Intervals}
\end{figure}
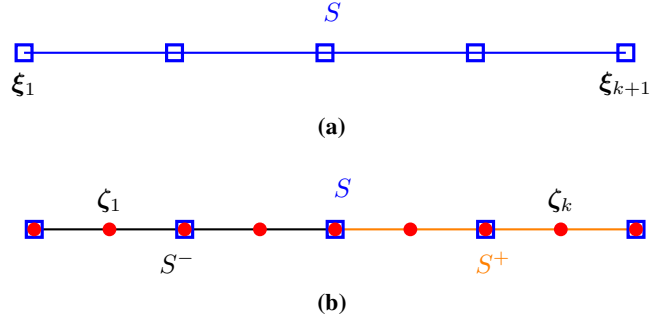

Next, let us clarify what happens when a hanging node is created. Let $S$ be an element edge that is being refined, i.e., split into two contiguous edges $S^-$ and $S^+$. Before the refinement, $S$ contains $k+1$ equally-spaced nodes $\bm{\xi}_n$, $n=1,\dots k+1$: the endpoints and the $k-1$ internal ones. After the refinement, $S$ contains $2k+1$ nodes, precisely $k+1$ equally-spaced nodes on each sub-edge $S^\pm$, with the midpoint in common; see Figure \ref{Figure:Intervals}. The spacing of the `old' nodes on $S$ was $\frac{|S|}k$ (where $|S|$ denotes the length of $S$), whereas the spacing of the `new' nodes is $\frac{|S|}{2k}$. Consequently,  $k+1$ of these nodes coincide with those initially on $S$, and the new nodes introduced in the refinement are only $k$. We will denote these latter by $\bm{\zeta}_i$, $i=1,\dots, k$.
This suggests the following definition.
\begin{defn}[closest neighbors of a node]\label{def:closest-nbr} With the previous notation, if $\bm{x} := \bm{\zeta}_i$ is created as the midpoint of the segment $[\bm{x}', \bm{x}'']:=[\bm{\xi}_{n_i},\bm{\xi}_{n_i+1}]$ for some $n_i$, we define the set $\mathcal{B}(\bm{x}):= \{ \bm{x}', \bm{x}''\}$.
\end{defn}

We are ready to give the announced definition of global index of a node of the triangulation $\mesh$.

\begin{figure}[t!]
\centering
\subfloat[]{
\scalebox{0.7}{
\begin{tikzpicture}
\draw [thick]plot coordinates{(0,0) (12*3/4,0)  (4*3/4,3) (0,0)};
\draw [thick]plot coordinates{(0,0) (12*3/4,0)  (4*3/4,-4.5*3/5) (0,0)};
\draw [thick]plot coordinates{ (6*3/4,0)  (4*3/4,-4.5*3/5) };
\draw [thick]plot coordinates{ (6*3/4,0)  (8*3/4,-2.25*3/5) };
\draw [thick]plot coordinates{ (9*3/4,0)  (8*3/4,-2.25*3/5) };
\draw[very thick, red] plot[mark =square, only marks, mark size=3pt] coordinates{(0,0)(3*3/4,0) (6*3/4,0)(9*3/4,0)(12*3/4,0)(4*3/4,-4.5*3/5)(2*3/4,-4.5*3/10)(8*3/4,-4.5*3/10)(5*3/4,-4.5*3/10)};
\draw[very thick, green] plot[mark = *, only marks, mark size=2.5pt] coordinates{(12*3/4,0)(9*3/4,0)(6*3/4,0) (4*3/4,-4.5*3/5)(5*3/4,-4.5*3/10)(7*3/4,-4.5*3/20)(10*3/4,-4.5*3/20)(8*3/4,-4.5*3/10)};
\draw[very thick, orange] plot[mark =diamond, only marks, mark size=4pt] coordinates{(6*3/4,0)(15/2*3/4,0)(9*3/4,0)(10.5*3/4,0)(12*3/4,0)(7*3/4,-4.5*3/20)(8.5*3/4,-4.5*3/20)(10*3/4,-4.5*3/20)(8*3/4,-4.5*3/10)};
\draw [very thick, blue] plot [mark=x, only marks, mark size=4pt] coordinates
{(0,0) (6*3/4,0) (12*3/4,0) (8*3/4,3/2)(2*3/4,3/2)  (4*3/4,3) (0,0) (2*3/4,-4.5*3/10)(4*3/4,-4.5*3/5)(8*3/4,-4.5*3/10)};
\coordinate [label= 0] (z1) at (0,0.1);
\coordinate [label= 0] (z1) at (6*3/4,0.1);
\coordinate [label= 0] (z1) at (12*3/4,0.1);
\coordinate [label= 1] (z1) at (3*3/4,0.1);
\coordinate [label= 1] (z1) at (9*3/4,0.1);
\coordinate [label= 2] (z1) at (7.5*3/4,0.1);
\coordinate [label= 2] (z1) at (10.5*3/4,0.1);
\end{tikzpicture}}}
\subfloat[]{
\scalebox{0.7}{
\begin{tikzpicture}
\draw [thick]plot coordinates{(0,0) (12*3/4,0)  (4*3/4,3) (0,0)};
\draw [thick]plot coordinates{(0,0) (12*3/4,0)  (4*3/4,-4.5*3/5) (0,0)};
\draw [thick]plot coordinates{ (6*3/4,0)  (4*3/4,-4.5*3/5) };
\draw [thick]plot coordinates{ (6*3/4,0)  (8*3/4,-2.25*3/5) };
\draw [thick]plot coordinates{ (9*3/4,0)  (8*3/4,-2.25*3/5) };
\draw[very thick, red] plot[mark =square, only marks, mark size=3pt] coordinates{(0,0)(2*3/4,0) (4*3/4,0)(6*3/4,0)(8*3/4,0)(10*3/4,0)(12*3/4,0)(1.333*3/4,-1.5*3/5)(2.667*3/4,-3*3/5)(4*3/4,-4.5*3/5)(6.667*3/4,-3*3/5)(9.333*3/4,-1.5*3/5)(4.667*3/4,-3*3/5)(5.333*3/4,-1.5*3/5)};
\draw[very thick, green] plot[mark = *, only marks, mark size=2.5pt] coordinates{(12*3/4,0)(10*3/4,0)(8*3/4,0)(6*3/4,0) (4*3/4,-4.5*3/5)(6.667*3/4,-3*3/5)(9.333*3/4,-1.5*3/5)(4.667*3/4,-3*3/5)(5.333*3/4,-1.5*3/5)(10.667*3/4,-0.75*3/5)(5.333*3/4,-3.75*3/5)(6.667*3/4,-0.75*3/5)(7.334*3/4,-1.5*3/5)(8*3/4,-2.25*3/5)};
\draw[very thick, orange] plot[mark =diamond, only marks, mark size=4pt] coordinates{(4.5,0)(6,0)(7*3/4,0)(8*3/4,0)(9*3/4,0)(10*3/4,0)(11*3/4,0)(12*3/4,0)(6.667*3/4,-0.75*3/5)(7.334*3/4,-1.5*3/5)(8*3/4,-2.25*3/5)(8.333*3/4,-1.5*3/5)(8.667*3/4,-0.75*3/5)(9.333*3/4,-1.5*3/5)(10.667*3/4,-0.75*3/5)};
\draw [very thick, blue] plot [mark=x, only marks, mark size=4pt] coordinates
{(0,0) (4*3/4,0)(8*3/4,0) (12*3/4,0) (9.333*3/4,1)(6.667*3/4,2)  (4*3/4,3)(2.667*3/4,2)(1.333*3/4,1) (0,0) (1.333*3/4,-1.5*3/5)(2.667*3/4,-3*3/5)(4*3/4,-4.5*3/5)(6.667*3/4,-3*3/5)(9.333*3/4,-1.5*3/5)};
\coordinate [label= 0] (z1) at (0,0.1);
\coordinate [label= 0] (z1) at (4*3/4,0.1);
\coordinate [label= 0] (z1) at (8*3/4,0.1);
\coordinate [label= 0] (z1) at (12*3/4,0.1);
\coordinate [label= 1] (z1) at (2*3/4,0.1);
\coordinate [label= 1] (z1) at (6*3/4,0.1);
\coordinate [label= 1] (z1) at (10*3/4,0.1);
\coordinate [label= 2] (z1) at (7*3/4,0.1);
\coordinate [label= 2] (z1) at (11*3/4,0.1);
\coordinate [label= 2] (z1) at (9*3/4,0.1);
\end{tikzpicture}}}
\caption{Triangulation after the three refinements in the case $k=2$ (a) and in the case $k=3$ (b). Blue crosses represent the original degrees of freedom. Red squares, green circles and orange triangles are used for the degrees of freedom of the first, second and third refinement, respectively. All nodes are proper, except those on the horizontal line, whose global index is reported.}
\label{Figure:Tri}
\end{figure}
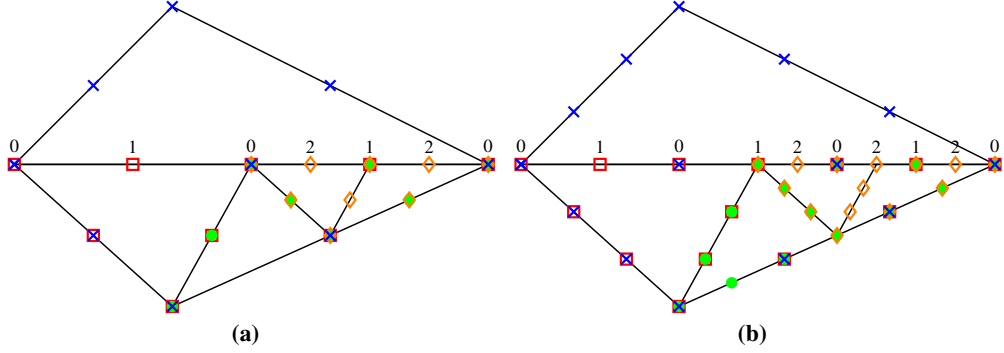

\begin{defn}[global index of a node] Given a node $\bm{x}\in\mathcal{N}$, we define its global index $\lambda$ recursively as follows:
\begin{itemize}
    \item if $\bm{x}$ is a proper node, then $\lambda(\bm{x}):=0$;
    \item if $\bm{x}$ is a hanging node, with $\bm{x}', \bm{x}'' \in \mathcal{B}(\bm{x})$, then set 
    \begin{linenomath}
    $$\lambda(\bm{x}) := \max\{\lambda(\bm{x}'), \lambda(\bm{x}'')\} +1.$$ \end{linenomath}
    \end{itemize}
\end{defn}
Figure \ref{Figure:Tri} shows the evolution of the global index after three refinements in the cases $k=2$ (a) and $k=3$ (b). We remark that, for instance, the midpoint of the horizontal edge is a proper node in case (a), and a hanging node in case (b).

The largest global index in $\mathcal{T}$ will be denoted by $\Lambda_\mathcal{T}:= \max_{\bm{x}\in \mathcal{N}}\{\lambda(\bm{x})\}$.
In this paper, as in \cite{AVEMstabfree},  we will consider sequences of successively refined triangulations $\{ \mesh \}$ whose global index does not blow up. 
\begin{assump}\label{AssumptionLambda}
There exists a constant $\Lambda>0$ such that, for any triangulation $\mesh$ generated by successive refinements of $\mesh_0$, it holds 
\begin{linenomath}
\begin{align*}
 \Lambda_\mathcal{T}\le \Lambda.
\end{align*}
\end{linenomath}
\end{assump}
Any such triangulation will be called {\em $\Lambda$-admissible}.
\section{Two key properties}\label{sec:TheKeyProperties}

In this section we discuss the validity of some results for the degree $k>1$ that will be used in the rest of the paper. We will highlight in particular the differences from the case $k=1$.

\begin{figure}
\centering
\begin{tikzpicture}[scale =0.9]
\draw [thick] plot coordinates{(8,0)  (3,2.5) (0,0)};
\draw [thick] plot coordinates{(8,0)  (4,-4*2/3) (0,0)};
\draw[ thick, red] plot coordinates{ (4,0)  (4,-4*2/3) };
\draw [thick, green]plot coordinates{ (4,0)  (2,-2*2/3) };
\draw [thick, green]plot coordinates{ (2,0)  (2,-2*2/3) };
\draw [thick, orange]plot coordinates{ (2,0)  (3,-1*2/3) };
\draw [ thick, orange]plot coordinates{ (3,0)  (3,-1*2/3) };
\draw [thick, blue]plot coordinates{ (0,0) (8,0)};

\draw[very thick, red] plot[mark =square, only marks, mark size=3pt] coordinates{(0,0)(4/3,0)(8/3,0)(4,0)(4 +4/3,0)(4 +8/3,0)(8,0)};
\draw[very thick, green] plot[mark =*, only marks, mark size=2.5pt] coordinates{(0,0)(4,0) (2,0)(2*1/3,0)(2*2/3,0)( 2+ 2/3, 0)(2+ 4/3,0)};
\draw[very thick, orange] plot[mark =diamond, only marks, mark size=4pt] coordinates{(2,0)(4,0) (2 +1/3,0) (2+2/3, 0) (3,0) (3+1/3,0) (3 +2/3,0) };
\draw [very thick, blue] plot [mark=x, only marks, mark size=4pt] coordinates {(0,0) (8/3,0)(16/3,0)(8,0)};

\coordinate [label= $|L|/3$] (z1) at (8/3,0.1);
\coordinate [label= $E$] (z1) at (3+1/3,1);
\coordinate [label= {[blue] $L$}] (z1) at (6,0.1);
\end{tikzpicture}
\caption{
The case $k=3$ with $3$ refinements of the edge $L$ (in blue) is shown. Red, green and orange lines are the lines needed to refine $L$ the first, the second and the third time respectively. Blue crosses are the degrees of freedom on $L$ of the function living on $E$. Red squares, green circles, orange diamonds are the degrees of freedom on $L$ generated after the first, the second and the third refinement of $L$.}
\label{Figure:Poincare2}
\end{figure}
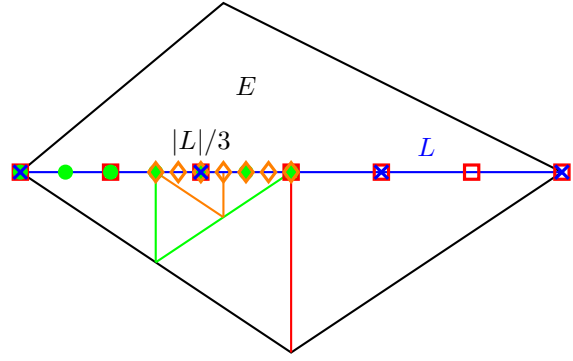

\begin{prop}[scaled Poincaré inequality in $\mathbb{V}_\mathcal{T}$]\label{ScaledPoincare}
    There exists a constant $C_P>0$, independent of $\mathcal{T}$, such that    
    \begin{linenomath}
    \begin{align}\label{eq:ScaledPoicare}
     &\sum_{E\in \mathcal{T}}h_E^2  \normLE{v}^2 \le C_P \normH{v}^2 &\forall v \in \mathbb{V}_\mathcal{T}\text{ such that } v(\bm{x})=0, \forall \bm{x}\in \mathcal{P}.
    \end{align}
    \end{linenomath}
\end{prop}    
\begin{proof}
Let $E\in \mathcal{T}$ be an element of the triangulation. 
    If $E$ is an element of the original partition $\mathcal{T}_0$, all its vertices are proper nodes. Otherwise, $E$ has been generated after some refinements by splitting an element $\widetilde{E}$ into two elements, $E$ and $E'$. Let $L$ be the common edge shared by $E$ and $E'$. If $L$ is not further refined, then all the nodes on $L$ are proper because they are shared by $E$ and $E'$. 
    If $L$ is refined and $k$ is even, then the midpoint of $L$ is a proper node.
    
    So, let us consider the case $k$ odd and let us assume that $L$ is refined $M\ge1$ times. We focus in particular on the internal node $\bar{\bm{x}}$ of $L$ is at distance $\frac{\abs{L}}{k}$ from one of the endpoints, Figure \ref{Figure:Poincare2} shows the case $k=3$. This point belongs to one of the $M+1$ intervals in which $L$ is refined, having width $\abs{L}/2^s$, for some $1\le s \le M$. We remark that $s$ depends on how $L$ has been refined (in the case of uniform refinements of $L$, one has $2^s= M +1$). We localize the chosen node $\bar{\bm{x}}$ in $L$ by defining an $m \ge 0$ such that 
 \begin{linenomath}
    \begin{align*}
        \frac{\abs{L}\;m}{2^s} \le \frac{\abs{L}}{k} \le \frac{\abs{L}(m +1)}{2^s},
    \end{align*}
 \end{linenomath}
    or, equivalently,
\begin{linenomath}    
    \begin{align}\label{ConditionPoincare}
        k\; m \le 2^s \le k\;(m+1).
    \end{align}
\end{linenomath}
    The interval going from $\frac{\abs{L}\;m}{2^s}$ to  $\frac{\abs{L}(m +1)}{2^s}$ is an edge for a smaller element $E'$, thus it contains $k-1$ internal nodes. Since they are equi-spaced, their positions are  at
\begin{linenomath}
    \begin{align*}
      &\frac{\abs{L}}{2^s} \left( m+\frac{n}{k}\right)  &\text{ with } n=0, \dots, k.
    \end{align*}
\end{linenomath}
    By taking $n = 2^s - m\;k$, which is compatible with conditions \eqref{ConditionPoincare}, we conclude that one of the internal nodes of $E'$ coincides with $\bar{\bm{x}}$.
    
    This guarantees that $E$ has at least one proper node $\bm{x}$ on its boundary. By hypothesis $v(\bm{x})=0$, and so we can apply the classical Poincar\'e inequality, 
    \begin{linenomath}
    \begin{align*}
        h^{-2}_E \normLE{v}^2 \lesssim \normHE{v}^2, 
    \end{align*}
    \end{linenomath}
    that concludes the proof.
    \end{proof}
    

\begin{remark}The previous proof exploits the fact that when $k>1$, each element of the triangulation contains at least a proper node. This differs from the case $k=1$ in which the edges do not contain internal nodes, and then elements with all hanging nodes as vertices are admissible.
As a further difference from the case $k=1$, we highlight that in Proposition \ref{ScaledPoincare} the constant $C_P$ does not depend on the constant $\Lambda$, whose existence has been introduced in Assumption \ref{AssumptionLambda}.
\rightline{$\square$}
\end{remark}

%

The next result we are going to establish is a hierarchical representation of the interpolation error $v - \IE v$ on the boundary $\partial E$ of an element $E \in \mesh$.
Assume that $v \in \VEk$, and let $L$ be an edge of $E$; for simplicity, in the sequel the restriction of $v$ to $L$, which is a piecewise polynomial of degree $k$, will be still denoted by $v$. The subsequent bisections of $L$ which generate the nodes in $\mathcal{N}_E \cap L$ allow us to write the difference $(v - \IE v)_{|L}$ telescopically as
\begin{linenomath}
  \begin{equation}\label{eq:telescoping}
(v - \IE v)_{|L} = \sum_{j=1}^{J_L} ({\cal I}_j - {\cal I}_{j-1})v \,; 
\end{equation}  
\end{linenomath}
here, ${\cal I}_0 = {\IE}_{|L}$, ${\cal I}_{J_L}$ is the identity operator, whereas ${\cal I}_jv$ for $1 \leq j \leq J_L-1$ is the piecewise polynomial of degree $k$ which interpolates $v$ on the partition of $L$ of level $j$, namely the partition formed by sub-edges of length $\leq \frac{|L|}{2^j}$.

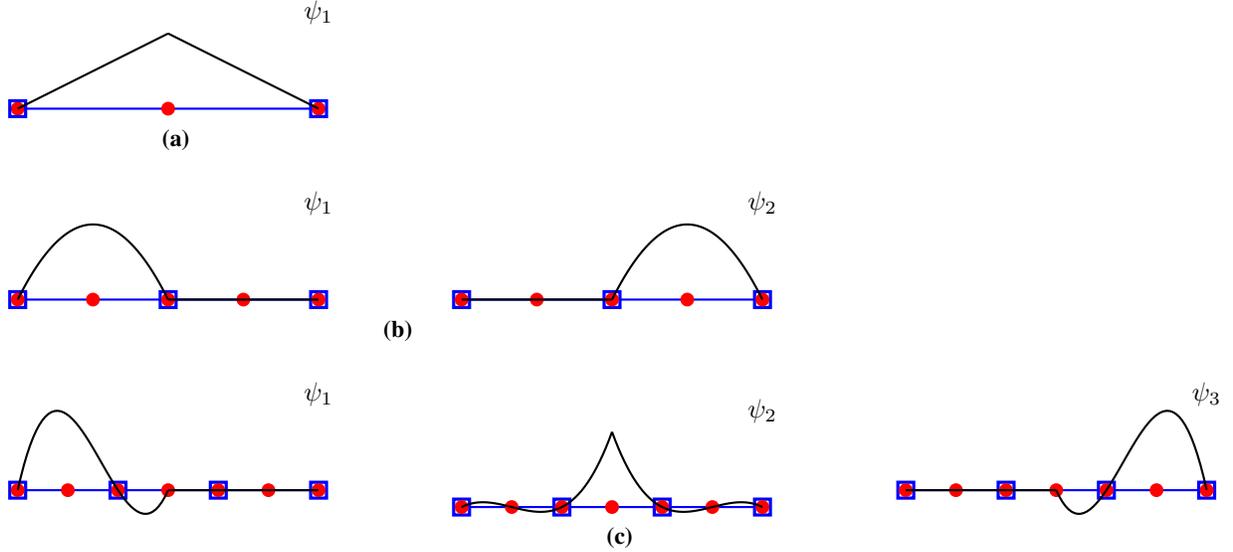
\begin{figure}[t!]

\subfloat[]{    
\begin{tikzpicture}
     \draw[thick,blue]plot coordinates {(0,0)(4,0)};
     \draw  [very thick, red] plot[mark =*, only marks, mark size=2pt]  coordinates{(0,0)(4,0)(2,0)};
    \draw [very thick,blue]plot [mark =square, only marks,   mark size=3pt] coordinates{(0,0)(4,0)};
    \draw [thick] plot [variable=\t, domain=0:2,samples=50]({\t}, {0.5*\t});
    \draw [thick] plot [variable=\t, domain=2:4,samples=50]({\t}, {-0.5*\t+2});
    \coordinate [label= $\psi_1$] (z1) at (4,1);
    \end{tikzpicture}
}\\
\subfloat[]{
\begin{tikzpicture}
     \draw[thick,blue]plot coordinates {(0,0)(4,0)};
     \draw [very thick,blue]plot [mark =square, only marks,   mark size=3pt] coordinates{(0,0)(2,0)(4,0)};
     \draw  [very thick, red] plot[mark =*, only marks, mark size=2pt]  coordinates{(0,0)(1,0)(2,0)(3,0)(4,0)};
    \draw [thick] plot [variable=\t, domain=0:2,samples=50]({\t}, {-1*\t*\t +2*\t});
    \draw [thick] plot [variable=\t, domain=2:4,samples=50]({\t}, 0);
    \coordinate [label= $\psi_1$] (z1) at (4,1);\end{tikzpicture}\hspace{1.4cm}
    \begin{tikzpicture}
     \draw[thick,blue]plot coordinates {(0,0)(4,0)};
     \draw [very thick,blue]plot [mark =square, only marks,   mark size=3pt] coordinates{(0,0)(2,0)(4,0)};
     \draw  [very thick, red] plot[mark =*, only marks, mark size=2pt]  coordinates{(0,0)(1,0)(2,0)(3,0)(4,0)};
    \draw [thick] plot [variable=\t, domain=0:2,samples=50]({\t}, 0 );
    \draw [thick] plot [variable=\t, domain=2:4,samples=50]({\t}, {-1*\t*\t +6*\t -8});
    
    \coordinate [label= $\psi_2$] (z1) at (4,1);
    \end{tikzpicture}
}\\
\subfloat[]{
\begin{tikzpicture}
     \draw[thick,blue]plot coordinates {(0,0)(4,0)};
     \draw [very thick,blue]plot [mark =square, only marks,   mark size=3pt] coordinates{(0,0)(4/3,0)(8/3,0)(4,0)};
     \draw  [very thick, red] plot[mark =*, only marks, mark size=2pt]  coordinates{(0,0)(2/3,0)(4/3,0)(6/3,0)(8/3,0)(10/3,0)(4,0)};
    \draw [thick] plot [variable=\t, domain=0:2,samples=50]({\t}, {27/16*(\t*(\t -4/3)*(\t-2))});
    \draw [thick] plot [variable=\t, domain=2:4,samples=50]({\t}, 0);
    
    \coordinate [label= $\psi_1$] (z1) at (4,1);\end{tikzpicture}\hspace{1.4cm}
    \begin{tikzpicture}
     \draw[thick,blue]plot coordinates {(0,0)(4,0)};
     \draw [very thick,blue]plot [mark =square, only marks,   mark size=3pt] coordinates{(0,0)(4/3,0)(8/3,0)(4,0)};
     \draw  [very thick, red] plot[mark =*, only marks, mark size=2pt]  coordinates{(0,0)(2/3,0)(4/3,0)(6/3,0)(8/3,0)(10/3,0)(4,0)};
    \draw [thick] plot [variable=\t, domain=0:2,samples=50]({\t}, {9/16*(\t*(\t -4/3)*(\t-2/3))} );
    \draw [thick] plot [variable=\t, domain=2:4,samples=50]({\t}, {-9/16*((\t-8/3)*(\t -4)*(\t-10/3))});
    
    \coordinate [label= $\psi_2$] (z1) at (4,1);\end{tikzpicture}\hspace{1.4cm}
    \begin{tikzpicture}
     \draw[thick,blue]plot coordinates {(0,0)(4,0)};
     \draw [very thick,blue]plot [mark =square, only marks,   mark size=3pt] coordinates{(0,0)(4/3,0)(8/3,0)(4,0)};
     \draw  [very thick, red] plot[mark =*, only marks, mark size=2pt]  coordinates{(0,0)(2/3,0)(4/3,0)(6/3,0)(8/3,0)(10/3,0)(4,0)};
    \draw [thick] plot [variable=\t, domain=0:2,samples=50]({\t}, 0 );
    \draw [thick] plot [variable=\t, domain=2:4,samples=50]({\t}, {-27/16*((\t-8/3)*(\t -4)*(\t-2))});
    
    \coordinate [label= $\psi_3$] (z1) at (4,1);
    \end{tikzpicture}
}
\caption{Blue square are the $k+1$ equi-spaced original nodes on the blue edge. Red points represent the nodes added after the refinement of the interval. Black lines show the shapes of the basis $\psi_i$, $i=1,\dots,k$, in the case $k=1$ (a), $k=2$ (b), $k=3$ (c). }
\label{Figure:BasisDefinition}
\end{figure}

\begin{figure}
\centering
\subfloat[]{
\begin{tikzpicture}[scale =0.9]
\draw [thick]plot coordinates{(0,0) (6,0)  (2,3) (0,0)};
\draw [thick]plot coordinates{(0,0) (6,0)  (2.5,-2) (0,0)};
\draw [thick]plot coordinates{(3,0)(2,3) };
\draw [very thick] plot [mark=*, only marks] coordinates
{(0,0) (3,0) (6,0) (4,1.5)  (2,3)(1,1.5)(2.5, -2)(2.5/2, -1) (2.5 + 3.5/2, -1)(2.5,1.5) };
\draw [very thick, red]plot [mark=*, only marks]coordinates{(1.5,0)(4.5,0)};
\coordinate [label= $\bm{\zeta}_1$] (z1) at (1.5,0.1);
\coordinate [label= $\bm{\zeta}_2$] (z1) at (4.5,0.1);
\coordinate [label= $\bm{\xi}_1$] (z1) at (0,-0.7);
\coordinate [label= $\bm{\xi}_2$] (z1) at (3.1,-0.7);
\coordinate [label= $\bm{\xi}_3$] (z1) at (6,-0.7);
\end{tikzpicture}}
\subfloat[]{
\begin{tikzpicture}[scale =0.9]
\draw [thick]plot coordinates{(0,0) (6,0)  (2,3) (0,0)};
\draw [thick]plot coordinates{(0,0) (6,0)  (2.5,-2) (0,0)};
\draw [thick]plot coordinates{(3,0)(2,3) };
\draw [very thick] plot [mark=*, only marks] coordinates
{(0,0) (2,0) (4,0) (6,0) (2/3,1) (2,3) (4/3,2) (2 +2/3,1) (2 + 1/3,2) (4+ 2/3,1) (3+1/3,2) (2.5, -2)(2.5/3, -2/3) (2.5 + 3.5*2/3 , -2/3) (2.5 + 3.5*1/3 ,-4/3) (2.5*2/3,-4/3) };
\draw [very thick, red]plot [mark=*, only marks]coordinates{(1,0)(3,0)(5,0)};
\coordinate [label= $\bm{\zeta}_2$] (z1) at (3,-0.7);
\coordinate [label= $\bm{\zeta}_1$] (z1) at (1,-0.7);
\coordinate [label= $\bm{\zeta}_3$] (z1) at (5,-0.7);
\coordinate [label= $\bm{\xi}_1$] (z1) at (0,0.1);
\coordinate [label= $\bm{\xi}_2$] (z1) at (2,0.1);
\coordinate [label= $\bm{\xi}_3$] (z1) at (4,0.1);
\coordinate [label= $\bm{\xi}_4$] (z1) at (6,0.1);
\end{tikzpicture}
}
\caption{Black points are the proper nodes. Red points represent the hanging nodes generated after a refinement. In (a) the case $k=2$ is showed, $\bm{\zeta}_1$ is the hanging node obtained after the refinement of $\bm{\xi}_1$ and $\bm{\xi}_3$ and it is the midpoint of $\bm{\xi}_1$ and $\bm{\xi}_2$. We notice that if we have called the other red point $\bm{\zeta}_2$, $\bm{\xi}_1$ and $\bm{\xi}_3$ would have been switched. Analogusly, (b) represents the case $k=3$. }
\label{Figure:LagrangeOperator}
\end{figure}
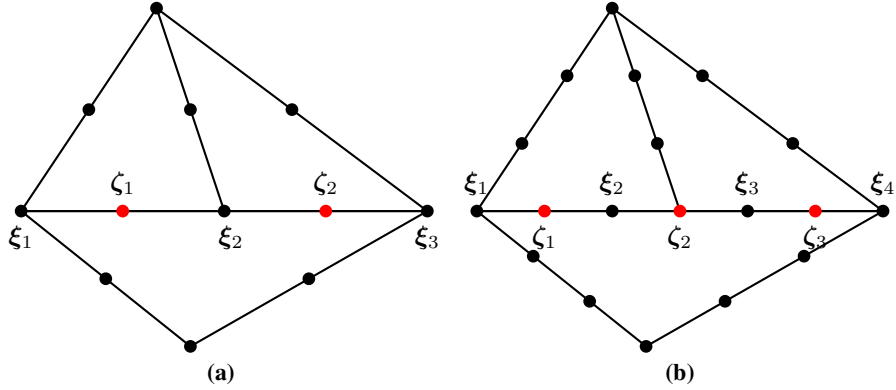

In order to understand the structure of the detail $({\cal I}_j - {\cal I}_{j-1})v$, assume that $S$ is a sub-edge of $L$ of length $= \frac{|L|}{2^{j-1}}$, which is split into two sub-edges $S^\pm$ of length $= \frac{|L|}{2^{j}}$ (see again Fig. \ref{Figure:Intervals}). On $S$ we have two interpolation operators, namely $\mathcal{I} := {{\cal I}_{j-1}}_{|S} 
:C^0(S)\rightarrow \mathbb{P}_k(S)$ and $\mathcal{I}^\vee := {{\cal I}_j}_{|L} :C^0(S)\rightarrow \mathbb{P}_{k}(S^-,S^+) = \left\{v\in C^0(S): v_{|S^-}\in \mathbb{P}_k(S^-) \text{ and } v_{|S^+}\in \mathbb{P}_k(S^+) \right\}$, which coincides with the interpolation operator $\mathcal{I}^-:C^0(S^-)\rightarrow \mathbb{P}_k(S^-)$ when restricted to $S^-$ and with the analogous operator $\mathcal{I}^+$ when restricted to $S^+$. 
With the notation introduced just before Definition \ref{def:closest-nbr}, we can quantify the discrepancy between the two interpolation operators by defining the  $k$ basis functions 
\begin{linenomath}
\begin{align*}
\psi_i \in \mathbb{P}_{k}(S^-,S^+) \text{ such that }
    \psi_i(\bm{x}) = \begin{cases}
         1 &\text{ if } \bm{x}= \bm{\zeta}_i,\\
         0 &\text{ if } \bm{x}= \bm{\zeta}_j, \;j\neq i,\\
         0 &\text{ if } \bm{x}= \bm{\xi}_n, \ n=1,\;\dots,k+1,
    \end{cases}
    \qquad 1 \leq i \leq k \,.
\end{align*}
\end{linenomath}
See Figure \ref{Figure:BasisDefinition} for a graphical representation of these functions in the cases $k=1$ (a), $k=2$ (b), $k=3$ (c).
Hence, the difference between the two interpolation operators on $S$ can be written as
\begin{linenomath}
\begin{align*}
    \mathcal{I}^\vee v  - \mathcal{I} v =\sum_{i=1}^k d(v,\bm{\zeta}_i)\psi_i,
\end{align*}
\end{linenomath}
where $d$ is defined as
\begin{linenomath}
\begin{align}\label{eq:interp_ric_3}
    d(v,\bm{\zeta}_i):= ( \mathcal{I}^\vee v - \mathcal{I}v)(\bm{\zeta}_i) = \left( v - \mathcal{I} v\right)(\bm{\zeta}_i).
\end{align}
\end{linenomath}
The values of $\mathcal{I}v$ at the $k$ nodes $\bm{\zeta}_i$ are a linear combination of the values of $\mathcal{I}v$ at the $k+1$ nodes $\bm{\zeta}_n$, where $\mathcal{I}v$ coincides with $v$. Thus, there exist coefficients $\alpha_{i,n}$ such that
\begin{linenomath}
    \begin{align}\label{interp_ric}
    & (\mathcal{I}v )(\bm{\zeta}_i)=\sum^{k+1}_{n=1} \alpha_{i,n} v(\bm{\xi}_n), &i=1,\dots,k.
\end{align}
\end{linenomath}
The explicit values of these coefficients in the case $k=2$ for the two new nodes $\bm{\zeta}_1$ and $\bm{\zeta}_2$ are given in these expressions:
\begin{linenomath}
\begin{align*}
&(\mathcal{I}v )(\bm{\zeta}_1) = \frac{3}{8}v(\bm{\xi}_1) + \frac{3}{4} v (\bm{\xi}_2) - \frac{1}{8} v(\bm{\xi}_3) \\
&(\mathcal{I}v )(\bm{\zeta}_2) =- \frac{1}{8} v(\bm{\xi}_1) + \frac{3}{4} v (\bm{\xi}_2) +\frac{3}{8}v(\bm{\xi}_3),    
\end{align*}
\end{linenomath}
where $\bm{\xi}_i \le \bm{\zeta}_i \le \bm{\xi}_{i+1}$, $i=1,2$. Similarly, in the case $k=3$, we get
\begin{linenomath}
\begin{align*}
&(\mathcal{I}v )(\bm{\zeta}_1) = \frac{5}{16}v(\bm{\xi}_1) + \frac{15}{16}v(\bm{\xi}_2) - \frac{5}{16}v(\bm{\xi}_3) + \frac{1}{16}v(\bm{\xi}_4),\\
&(\mathcal{I}v )(\bm{\zeta}_2) = -\frac{1}{16}v(\bm{\xi}_1) + \frac{9}{16}v(\bm{\xi}_2) + \frac{9}{16}v(\bm{\xi}_3) - \frac{1}{16}v(\bm{\xi}_4),\\
&(\mathcal{I}v )(\bm{\zeta}_3) = \frac{1}{16}v(\bm{\xi}_1) - \frac{5}{16}v(\bm{\xi}_2) + \frac{15}{16}v(\bm{\xi}_3) + \frac{5}{16}v(\bm{\xi}_4),
\end{align*}
\end{linenomath}
where again $\bm{\xi}_i \le \bm{\zeta}_i \le \bm{\xi}_{i+1}$, $i=1,2,3$. Figure \ref{Figure:LagrangeOperator} shows both cases.  We notice that the coefficients $\alpha_{i,n}$ depend only on the relative positions of the nodes on $S$, not on the level $j$ of refinement. 

Summarizing, at the level $j$ of refinement of the edge $L$, we get
\begin{linenomath}
$$
    (\mathcal{I}_j  - \mathcal{I}_{j-1}) v =\sum_{\bm{x} \in \mathcal{H}_{L,j}} d(v,\bm{x})\psi_{\bm{x}}\,,
$$
\end{linenomath}
where $\mathcal{H}_{L,j}$ is the set of hanging nodes on $L$ created at the level $j$ of refinement, whereas
\begin{linenomath}
$$
d(v,\bm{x}) =  (\mathcal{I}_j v - \mathcal{I}_{j-1}v)(\bm{x}) = \left( v - \mathcal{I}_{j-1} v\right)(\bm{x}).
$$\,
\end{linenomath}
Summing-up over the levels and recalling \eqref{eq:telescoping}, we obtain
\begin{linenomath}
$$
(v - \IE v)_{|L} = \sum_{\bm{x} \in \mathcal{H}_{L}} d(v,\bm{x})\psi_{\bm{x}} \,.
$$
\end{linenomath}
where $ \mathcal{H}_{L} =  \mathcal{H}_{E} \cap L$, whence
\begin{linenomath}
$$
(v - \IE v)_{|\partial E} = \sum_{\bm{x} \in \mathcal{H}_{E}} d(v,\bm{x})\psi_{\bm{x}}\,.
$$\end{linenomath}
We now introduce the subspace of $\VEk$
\begin{linenomath}
\begin{align*}
    X_E:= \left\{w \in \VEk: w(\bm{x})=0\; \forall \bm{x} \in \mathcal{P}_E, \text{ and } \bm{\mu}_p(w,E)=\bm{0}, \; 0\le p\le k-3 \right\},
\end{align*}
\end{linenomath}
which contains $v-\IE v$ by definition of $\IE$. On $X_E$, we have two norms, namely the seminorm $|w |_{1,E}$ (which is a norm on $X_E$ due to the vanishing of $w$ at the three vertices of $E$) and the norm
\begin{linenomath}
\begin{align*}
    &[\![w]\!]_{X_E}:= \left(\sum_{\bm{x} \in \mathcal{H}_E} d^2 (w,\bm{x}) +| \bm{\mu}_{k-2}(E,w)|^2 \right)^{1/2} \,.
\end{align*}
\end{linenomath}
Note that, due to Assumption \ref{AssumptionLambda}, the dimension of $X_E$ is uniformly bounded by a constant depending on $\Lambda$; furthermore, the number of possible patterns of hanging nodes on $\partial E$, which determine the details $d(w,\bm{x})$, is also bounded in terms of $\Lambda$. As a consequence, the two norms are equivalent, with equivalence constants depending on $\Lambda$. Therefore,
\begin{linenomath}
$$
\sum_{\bm{x} \in \mathcal{H}_E} d^2 (w,\bm{x}) \leq [\![w]\!]_{X_E}^2 \simeq |w |_{1,E}^2 \qquad \forall w \in X_E\,.
$$
\end{linenomath}
Since $v - \IE v \in X_E$ and $d (v - \IE v,\bm{x}) = d(v, \bm{x})$ for any $\bm{x} \in \mathcal{H}_E$, we obtain
\begin{linenomath}
$$
\sum_{\bm{x} \in \mathcal{H}_E} d^2 (v,\bm{x}) \lesssim |v-\IE v |_{1,E}^2 \,.
$$
\end{linenomath}
Summing-up over all the elements of the triangulation, we arrive at the following result.

\begin{lemma}[global interpolation error vs hierarchical errors]\label{globInterVsHierErr} 
There exists a constant $C_D > 0$ depending on $\Lambda$ but independent of the triangulation $\mathcal{T}$ such that
\begin{linenomath}
\begin{equation}\label{eq:SimDelta}
     \sum_{\bm{x} \in \mathcal{H}}d^2(v,\bm{x})\le C_D \normHT{v - \mathcal{I}_\mathcal{T}v}^2  \qquad  \forall v \in\mathbb{V}_\mathcal{T} \,.
\end{equation}
\end{linenomath}
\end{lemma}

Next, we introduce the interpolation operator
\begin{linenomath}
\begin{equation}\label{eq:Imesh0}
\Imeshz : \Vmesh \to \Vmeshz \,,
\end{equation}
\end{linenomath}
where $\Vmeshz$ is defined in \eqref{def:V0T}, by the following conditions:
\begin{itemize}
\item $(\Imeshz v)(\bm{x}) = v(\bm{x})$ for all $\bm{x} \in {\cal P}$\,,
\item $\bm{\mu}_p(E,\Imeshz v) =  \bm{\mu}_p(E, v)$ for all $0 \leq p \leq k-3$ and for all $E \in \mesh$.
\end{itemize}
These conditions uniquely identify $\Imeshz v$. Indeed, if $\bm{x} \in {\cal H}$ is generated by a refinement of level $j$ of an edge $L$ (say, $\bm{x}=\bm{\zeta}_i$ with the notation introduced before Definition \ref{def:closest-nbr}), then $(\Imeshz v)(\bm{x}) $ can be expressed in terms of the values of $\Imeshz v$ at the $k+1$ nodes (say, $\bm{\xi}_n$) created at the previous levels of refinement of $L$, using the same coefficients as in formula \eqref{interp_ric}, i.e.,
\begin{linenomath}
\begin{align}\label{interp_ric_2}
    & (\Imeshz v)(\bm{\zeta}_i)=\sum^{k+1}_{n=1} \alpha_{i,n} (\Imeshz v)(\bm{\xi}_n),  \quad i=1,\dots,k\,; 
\end{align}
\end{linenomath}
and so on recursively. 

The following result provides a representation of the error $\Imesh v - \Imeshz v$.
\begin{lemma}\label{lemma:iterp_err}
It holds
\begin{linenomath}
    \begin{align*}
    \normHT{\mathcal{I}_\mathcal{T} v -\mathcal{I}^0_\mathcal{T} v }^2 \simeq\sum_{\bm{x}\in \mathcal{H}}\delta^2(v,\bm{x}) \qquad \forall v \in \Vmesh \,,
\end{align*}
\end{linenomath}
where $\delta(v,\bm{x}):= v(\bm{x}) - (\Imeshz v)(\bm{x})$.
\end{lemma}
\begin{proof}  
Consider an element $E \in \mesh$.
Recall that by construction  it holds $\bm{\mu}_p(E,\mathcal{I}_E v) = \bm{\mu}_p(E,v) = \bm{\mu}_p(E,\Imeshz v)$, whence
$\bm{\mu}_p(\mathcal{I}_E v - \Imeshz v,E)=\bm{0}$ for all $0 \leq p \leq k-3$. Consequently,
\begin{linenomath}
\begin{align*}
    \normHE{\mathcal{I}_E v - \Imeshz v}^2 \simeq \sum_{\bm{x}\in \mathcal{P}_E} \abs{\left( \mathcal{I}_E v - \Imeshz v\right)(\bm{x})}^2.
\end{align*}
\end{linenomath}
If $\bm{x}\in\mathcal{P}_E$, $(\mathcal{I}_E v)(\bm{x})=v(\bm{x})$, hence
\begin{linenomath}
\begin{align*}
    \normHE{\mathcal{I}_E v - \Imeshz v}^2 \simeq \sum_{\bm{x}\in \mathcal{P}_E} \abs{\left(  v - \Imeshz v\right)(\bm{x})}^2.
\end{align*}
\end{linenomath}
Summing on all the elements of the partition, we get
\begin{linenomath}
\begin{align*}
    & \sum_{E\in \mathcal{T}}\normHE{\mathcal{I}_E v - \Imeshz v}^2 \simeq \sum_{\bm{x}\in \mathcal{N}}\abs{\left(  v - \Imeshz v\right)(\bm{x})}^2 \simeq \sum_{\bm{x}\in \mathcal{H}}\abs{\left(  v - \Imeshz v\right)(\bm{x})}^2,
\end{align*}
\end{linenomath}
since if $\bm{x}\in\mathcal{P}$, $(\mathcal{I}^0_E v)(\bm{x})=v(\bm{x})$. This concludes the proof. 
\end{proof}

Concatenating Lemma \ref{globInterVsHierErr} and Lemma \ref{lemma:iterp_err}, we can prove the second key property of this section.
\begin{prop}[comparison between interpolation operators]\label{Prop:InterpolationOperators}
Let $\mathcal{T}$be $\Lambda$-admissible. Then, there exists a constant $C_I>0$, depending on $\Lambda$, but independent of $\mathcal{T}$, such that
\begin{linenomath}
\begin{align*}
    &\normH{v - \mathcal{I}^0_\mathcal{T} v} \le C_I \normHT{v- \mathcal{I}_\mathcal{T} v}  \qquad \forall v \in \mathbb{V}_\mathcal{T}.
\end{align*}
\end{linenomath}
\end{prop}
\begin{proof} Given a function $v \in \mathbb{V}_ {\mathcal{T}}$, by the triangle inequality
\begin{linenomath}
\begin{align*}
    \normH{v-\mathcal{I}^0_\mathcal{T}v}= |v-\mathcal{I}^0_\mathcal{T}v|_{1,\mathcal{T}}\le |v-\mathcal{I}_\mathcal{T}v|_{1,\mathcal{T}} +|\mathcal{I}_\mathcal{T}v-\mathcal{I}^0_\mathcal{T}v|_{1,\mathcal{T}} \,,
\end{align*}
\end{linenomath}
so it is enough to bound the last norm on the right-hand side. To this end, considering the vectors
\begin{linenomath}
\begin{align*}
    &\bm{\delta} = \left(\delta(\bm{x})\right)_{\bm{x}\in \mathcal{H}} :=  \left(\delta(v,\bm{x})\right)_{\bm{x}\in \mathcal{H}}\,,   \qquad  \bm{d} = \left(d(\bm{x})\right)_{\bm{x}\in \mathcal{H}} := \left(d(v,\bm{x})\right)_{\bm{x}\in \mathcal{H}},
\end{align*}
\end{linenomath}
and recalling the two Lemmas, the proof can be concluded if we show that
\begin{linenomath}
\begin{align*}
    \norm{\bm{\delta}}_{l^2(\mathcal{H})}\lesssim \norm{\bm{d}}_{l^2(\mathcal{H})}.
\end{align*}
\end{linenomath}

Given $\bm{x} \in {\cal H}$, assume that it is generated by a refinement of level $j$ of an edge $L$ (say, $\bm{x}=\bm{\zeta}_i$ with the notation introduced before Definition \ref{def:closest-nbr}). Writing $v^*:=\mathcal{I}^0_\mathcal{T}v$ for short, and exploiting formulas \eqref{eq:interp_ric_3} and \eqref{interp_ric}, we get
\begin{linenomath}
\begin{equation}\label{eq:delta-d}
\begin{split}
    \delta(\bm{\zeta}_{i})&= v(\bm{\zeta}_{i})- v^*(\bm{\zeta}_{i})= v(\bm{\zeta}_{i}) - \sum^{k+1}_{n=1} \alpha_{i,n} v^*(\bm{\xi}_n) \\&= v(\bm{\zeta}_{i}) - \sum^{k+1}_{n=1} \alpha_{i,n} v(\bm{\xi}_n)- \sum^{k+1}_{n=1} \alpha_{i,n}\left( v^*(\bm{\xi}_n) - v(\bm{\xi}_n) )\right)
    \\&= d(\bm{\zeta}_{i}) + \sum^{k+1}_{n=1} \alpha_{i,n} \delta(\bm{\xi}_n)\,.
\end{split}
\end{equation}
\end{linenomath}
Thus, we can build a matrix $\bm{W}: l^2(\mathcal{H})\rightarrow l^2(\mathcal{H})$ such that $\bm{\delta}=\bm{W}\bm{d}$, and we just need to prove that
\begin{linenomath}
\begin{align*}
    ||\bm{W}||_2\lesssim 1.
\end{align*}
\end{linenomath}
We now organize the hanging nodes with respect to the global index $\lambda \in [1,\Lambda_\mathcal{T}]$. Calling $\mathcal{H}_\lambda = \{\bm{x} \in \mathcal{H}: \lambda(\bm{x})=\lambda\}$, and $\mathcal{H} =\bigcup_{1\le \lambda\le \Lambda_\mathcal{T}}\mathcal{H}_\lambda$, the matrix $\bm{W}$ can be factorized in lower triangular matrices $\bm{W}_\lambda$, that change the nodes of level $\lambda$, leaving the others unchanged. In particular,
\begin{linenomath}
\begin{align*}
    \bm{W} = \bm{W}_{\Lambda_\mathcal{T}} \bm{W}_{\Lambda_{\mathcal{T}}-1}...  \bm{W}_{2} \bm{W}_{1},
\end{align*}
\end{linenomath}
where $\bm{W}_{1}$ is just the identity matrix $\bm{I}$, whereas each other matrix $\bm{W}_\lambda$ differs from the identity only in  the rows of block $\lambda$. In each of these rows, all entries are zero,  but the entries  $\alpha_{i,n}$ in the off-diagonal part and $1$ on the diagonal.
In order to estimate $\bm{W}_\lambda$, we use the H\"older inequality $||\bm{W}_\lambda||^2_2\le ||\bm{W}_\lambda||_1||\bm{W}_\lambda||_\infty$.
From the construction of $\bm{W}_\lambda$ have that 
\begin{linenomath}
\begin{align*}
 &||\bm{W}_\lambda||_{\infty}\le \max_n\left \{ \sum_{i=1}^{k+1} \abs{\alpha_{i,n}} \right \}+1=: \beta_1\,, &||\bm{W}_\lambda||_{1}\le 5 \;k\; \max_{i,n}{\abs{\alpha_{i,n}}}+1=: \beta_2 \,,
\end{align*}
\end{linenomath}
where in the last inequality it has been used the fact that a hanging node of global index $<\lambda$ may appear at most 5 times on the right-hand side of  \eqref{eq:delta-d}, since at most five edges meet at a node \cite[Proposition 3.2]{AVEMstabfree}.
These bring us to the following bound
\begin{linenomath}
\begin{equation*}
    ||\bm{W}||_2\le \prod_{2\le \lambda\le \Lambda_\mathcal{T}}||\bm{W}_\lambda||_2 \le\left(\beta_1 \cdot \beta_2\right)^{\frac{\Lambda-1}{2}}.
\end{equation*}
\end{linenomath}
and the proof is concluded.
\end{proof}

 \section{A posteriori error estimator}\label{sec:aPosterioriErrorEstimator}
  
With the aim of discussing the a posteriori error analysis, and following \cite{Cangiani}, we define the a posteriori error estimators, starting from the internal residual over an element $E$, i.e.,
\begin{linenomath}
\begin{align}\label{ResLoc}
    r_\mathcal{T}(E; v, \mathcal{D}):= f_E + \nabla \cdot \left(A_E \Pi^0_{k-1} \nabla v \right) - c_E \Pi^0_k v,
\end{align}
\end{linenomath}
for any $v\in \mathbb{V}_{E,k}$. We highlight that in the case $k=1$, with piecewise constant data, the diffusion term in the residual vanishes.
Furthermore, we define the jump residual over $e$, where $e$ is an edge shared by two elements $E_1$ and $E_2$ of the partition $\mathcal{T}$, as
\begin{linenomath}
\begin{align*}
    j_\mathcal{T}(e;v, \mathcal{T}):= [[ A  \Pi^0_{k-1} \nabla v]]_e = (A_{E_1} \Pi^0_{k-1} \nabla v|_{E_1})\cdot \bm{n_1}+( A_{E_2} \Pi^0_{k-1} \nabla v|_{E_2})\cdot \bm{n_2},
\end{align*}
\end{linenomath}
where $\bm{n_i}$ denotes the unit normal vector to $e$ pointing outward with respect to $E_i$; we set $j_\mathcal{T}(e;v)=0$ of $e \in \partial \Omega$.
Then, let the local residual estimator associated with $E$ be
\begin{linenomath}
\begin{align}\label{etaLocT}
    \eta^2_\mathcal{T}(E; v, \mathcal{D}):= h^2_E ||r_\mathcal{T}(E; v, \mathcal{D})||^2_{0,E} + \frac{1}{2} \sum_{e \in \mathcal{E}_E}h_E ||j_\mathcal{T}(e; v, \mathcal{D})||^2_{0,e},
\end{align}
\end{linenomath}
and the global residual estimator as the sum of the local residuals
\begin{linenomath}
\begin{align*}
    \eta^2_\mathcal{T}(v, \mathcal{D}):= \sum_{E \in \mathcal{T}} \eta^2_{\mathcal{T}}(E;v, \mathcal{D}).
\end{align*}
\end{linenomath}
In contrast to what has been done for the case $k=1$, we also need to introduce the \emph{virtual inconsistency terms}, defined by
\begin{linenomath}
\begin{equation}\label{eq:def-inconsistency}
\begin{split}
    \Psi^2_{\mesh,A}(E; v, \mathcal{D})&:=||(I-\Pi^0_{k-1})
    (A_E\Pi^0_{k-1}\nabla v) 
    ||^2_{0,E}, \\
    \Psi^2_{\mesh,c} (E; v, \mathcal{D})&:=h^2_E ||\left (I-\Pi^0_k)(c_E \Pi^0_k v \right)||^2_{0,E},
    \end{split}
\end{equation}
\end{linenomath}
as well as their sum
\begin{linenomath}
    \begin{align}\label{def:PsiAc}
    &\Psi^2_\mathcal{T}(v, \mathcal{D}):= \sum_{E\in \mathcal{T}}\Psi_\mesh^2(E; v, \mathcal{D})  := \sum_{E\in \mathcal{T}} \Psi^2_{\mesh,A}(E; v, \mathcal{D}) + \Psi^2_{\mesh, c}(E; v, \mathcal{D}).
\end{align}
\end{linenomath}

 \section{A posteriori error estimates}\label{sec:aPosterioriErrorEstimates}
 
In this section we present one of the main results of this paper, a stabilization-free a posteriori error bound. In this view, we firstly start by introducing the classical Clément operator upon the space $\mathbb{V}^0_\mesh$, $\Tilde{\mathcal{I}}^0_\mathcal{T}: \mathbb{V}\rightarrow \mathbb{V}^0_\mesh$; it is defined at the proper nodes on the skeleton of $\mesh$ as the average of the target function on the support of the associated basis functions, whereas the internal moments (if any) coincide with those of the target function.

The scaled Poincaré inequality (Proposition \ref{ScaledPoincare}) and Proposition \ref{Prop:InterpolationOperators} guarantee the validity of the error estimate for $\Tilde{\mathcal{I}}^0_\mathcal{T}$. Given these propositions, its proof does not involve the polynomial degree $k$, hence, it does not change with respect  to the one presented in \cite{AVEMstabfree}.
 \begin{lemma}[Clément interpolation estimate]\label{Lemma:ClementEstimate}
     $\forall v \in \mathbb{V}$, it holds
\begin{linenomath}
\begin{align*}
    \sum_{E\in \mathcal{T}}h^{-2}_E \normLE{v - \Tilde{\mathcal{I}}^0_\mathcal{T}v}^2\lesssim \normH{v}^2,
\end{align*} 
\end{linenomath}
where the hidden constant depends on $\Lambda$ but not on $\mathcal{T}$.
 \end{lemma}
 
We can now prove the following results, which is similar to Theorem 13 in \cite{Cangiani}, but with a slightly modified proof.

\begin{prop}[upper bound]\label{upperBound}
    There exists a constant $C_\text{apost}>0$, independent of $u$, $\mathcal{T}$, $u_\mathcal{T}$ and $\gamma$, such that
\begin{linenomath}
\begin{align}
    |u -u_\mathcal{T}|^2_{1,\Omega} \le C_\text{apost} \left(\eta^2_\mathcal{T}(u_\mathcal{T}, \mathcal{D})+ S_\mathcal{T}(u_\mathcal{T}, u_\mathcal{T})\right).
\end{align}
\end{linenomath}    
\begin{proof}
For any $v \in\mathbb{V}$, using the definition of Problem \eqref{Variational_Problem}, we have that
\begin{linenomath}
\begin{align*}
 \mathcal{B}(u - u_\mathcal{T}, v)&=\mathcal{B}(u,v) - \mathcal{B}(u_\mathcal{T},v)- (f,v_\mathcal{T})+ \mathcal{B}(u,v_\mathcal{T})\\&=(f, v-v_\mathcal{T})-\mathcal{B}(u_\mathcal{T},v)+\mathcal{B}(u,v_\mathcal{T})-\mathcal{B}(u_\mathcal{T},v_\mathcal{T})+\mathcal{B}(u_\mathcal{T},v_\mathcal{T})\\
 &=\left((f, v-v_\mathcal{T})-\mathcal{B}(u_\mathcal{T},v-v_\mathcal{T})\right)+\mathcal{B}(u-u_\mathcal{T},v_\mathcal{T})=: I +II,
\end{align*}
\end{linenomath}
where $v_\mathcal{T}:= \Tilde{\mathcal{I}}^0_\mathcal{T} v \in \Vmesh^0$.
The first term can be written as
\begin{linenomath}
\begin{align*}
    I&=
    \sum_{E\in \mathcal{T}}\left\{\int_E f_E(v-v_\mathcal{T})- \int_E A_E \nabla u_\mathcal{T} \cdot \nabla (v-v_\mathcal{T})- \int_E c_E u_\mathcal{T} (v-v_\mathcal{T})\right\}\\
    &= \sum_{E\in \mathcal{T}}\left\{\int_E f_E(v-v_\mathcal{T})- \int_E  \left(A_E \Pi^0_{k-1}\nabla u_\mathcal{T}\right) \cdot\nabla (v-v_\mathcal{T})- \int_E \left(c_E \Pi^0_{k}u_\mathcal{T}\right) (v-v_\mathcal{T})\right\}\\
     & \qquad +\sum_{E\in \mathcal{T}}\left\{\int_E \left( A_E(\Pi^0_{k-1} -I)\nabla u_\mathcal{T} \right )\cdot \nabla(v-v_\mathcal{T})+ \int_E \left(c_E( \Pi^0_k -I )u_\mathcal{T} \right)(v-v_\mathcal{T})\right\} = : I_1 +I_2.
\end{align*}
\end{linenomath}
The addend $I_1$ can be expressed as 
\begin{linenomath}
    \begin{align*}
    I_1 &=\sum_{E\in \mathcal{T}}\left\{\int_E \left ( f_E + \nabla \cdot \left( A_E \Pi^0_{k-1}\nabla u_\mathcal{T} \right)-  c_E \Pi^0_{k}u_\mathcal{T}\right) (v-v_\mathcal{T}) \right\}\\
    & \qquad +  \sum_{E\in \mathcal{T}} \int_{\partial E}\bm{n} \cdot \left( A_E \Pi^0_{k-1} 
 \nabla u{_\mathcal{T}}\right) (v-v_\mathcal{T}), 
 \end{align*}
\end{linenomath}
which can be bounded by using Lemma \ref{Lemma:ClementEstimate},
\begin{linenomath}
\begin{align*}
    |I_1|\lesssim \eta_\mathcal{T} (u_\mathcal{T}, \mathcal{D})|v|_{1,\Omega}.
\end{align*}
\end{linenomath}
On the other hand, noting that 
\begin{linenomath}
\begin{align}\label{property:Pi0}
\normLE{(I - \Pi^0_{k-1}) \nabla u_\mathcal{T}} =\normLE{(I - \Pi^0_{k-1})\nabla (I- \Pi^0_k)u_\mathcal{T}}\le \normLE{\nabla (I- \Pi^0_k)u_\mathcal{T}}  
\end{align}
\end{linenomath}
and applying again Lemma \ref{Lemma:ClementEstimate}, the addend $I_2$ can be bounded as follows:
\begin{linenomath}
\begin{align*}
    |I_2| &\le \left(\sum_{E\in \mathcal{T}}h^2_E|| \nabla\cdot\left( A_E(I-\Pi^0_{k-1})\nabla u_\mathcal{T} \right )||^2_{0,E}+ h^2_E ||\left(c_E (I-\Pi^0_k)u_\mathcal{T} \right)||^2_{0,E}\right)^{1/2}\\
    & \qquad \times \left(\sum_{E\in \mathcal{T}}h^{-2}_E||v - v_\mathcal{T}||^2_{0,E}\right)^{1/2}
    \\&\lesssim \left(\sum_{E\in \mathcal{T}}|| A_E(I-\Pi^0_{k-1})\nabla u_\mathcal{T} ||^2_{0,E}+ h^2_E ||c_E(I-\Pi^0_k)u_\mathcal{T} ||^2_{0,E}\right)^{1/2} |v|_{1,\Omega} \\
    &\lesssim \left(\sum_{E\in \mathcal{T}}|| \nabla  (u_\mathcal{T}-\Pi^0_{k}u_\mathcal{T} )||^2_{0,E}+ h^2_E ||( u_\mathcal{T}-\Pi^0_ku_\mathcal{T} ) ||^2_{0,E}\right)^{1/2} |v|_{1,\Omega} 
    \\
    &\lesssim \left(S_\mathcal{T}(u_\mathcal{T}, u_\mathcal{T})\right)^{1/2} |v|_{1,\Omega} \,.
\end{align*}
\end{linenomath}
Looking now at the term $II$, we have by Lemma \ref{Lemma:QuasiOrtog}
\begin{linenomath}
    \begin{equation*}
        | \mathcal{B}(u - u_\mathcal{T}, v)| \lesssim S_\mathcal{T}(u_\mesh, u_\mesh)^{1/2} |v|_{1,\Omega} \,.
    \end{equation*}
\end{linenomath}
Finally, by taking $v:= u- u_\mathcal{T} \in \mathbb{V}$, we get
\begin{linenomath}
\begin{align*}
    \mathcal{B}(u - u_\mathcal{T}, u - u_\mathcal{T})\lesssim \left(\eta_\mathcal{T} (u_\mathcal{T}, \mathcal{D})  + S_\mathcal{T}(u_\mathcal{T}, u_\mathcal{T})^{1/2} \right)|u - u_\mathcal{T}|_{1,\Omega},
\end{align*}
\end{linenomath}
which, using the coercivity of $\mathcal{B}$,  concludes the proof.
\end{proof}
\end{prop}
We now report a bound for the local residual estimator, proved in \cite{Cangiani}(Theorem 16).

\begin{prop}[local lower bound] There holds
\begin{linenomath}
\begin{align*}
    \eta^2_\mathcal{T}(E;u_\mathcal{T}, \mathcal{D})\lesssim \sum_{E'\in w_E}\left(|u-u_\mathcal{T}|^2_{1,E'} + S_{E'}(u_\mathcal{T},u_\mathcal{T})\right)
\end{align*}
\end{linenomath}
    where $w_E:=\{E':|\partial E \cap \partial E'| \neq 0\}$. The hidden constant is independent of $\gamma$, $h$, $u$ and $u_\mathcal{T}$. 
\end{prop}
Summing on all the elements of the partition, we get the following corollary. 
\begin{cor}[global lower bound]\label{cor:globalLoweBound} There exists a constant $c_\text{apost}>0$, independent of $u$, $\mathcal{T}$, $u_\mathcal{T}$ and $\gamma$, such that
\begin{linenomath}
\begin{align*}
    c_\text{apost} \;\eta^2_\mathcal{T}(u_\mathcal{T}, \mathcal{D})\le |u - u_\mathcal{T}|^2_{1,\Omega} + S_\mathcal{T}(u_\mathcal{T},u_\mathcal{T}).
\end{align*}
\end{linenomath}
\end{cor}

In the following proposition we present a bound of the stabilization term. We remark  that in the case $k=1$ the inconsistency term does not appear.

\begin{prop}[bound of the stabilization term]\label{bSt} There exists a constant $C_B > 0$ independent of $\mathcal{T}$, $u_\mathcal{T}$ and $\gamma$, such that
\begin{linenomath}
\begin{align}
     \gamma^2 S_\mathcal{T}(u_\mathcal{T},u_\mathcal{T})\le C_B  \left(\eta^2_\mathcal{T}(u_\mathcal{T}, \mathcal{D}) + \Psi^2_\mathcal{T}(u_\mathcal{T}, \mathcal{D})\right).
\end{align}
\end{linenomath}
  \end{prop}
\begin{proof}
From the definition \eqref{DefBT} of the form $\mathcal{B}_\mathcal{T}$ and from \eqref{Discrete_Variazional_Problem}, $\forall w \in \mathbb{V}^0_\mathcal{T}$ it holds 
\begin{linenomath}
\begin{align*}
    \gamma S_\mathcal{T}(u_\mathcal{T}, u_\mathcal{T})&= \gamma S_\mathcal{T}(u_\mathcal{T}, u_\mathcal{T}-w)\\&= \mathcal{B}_\mathcal{T}(u_\mathcal{T},u_\mathcal{T}-w) - a_\mathcal{T}(u_\mathcal{T}, u_\mathcal{T}-w) - m_\mathcal{T}(u_\mathcal{T},u_\mathcal{T}-w)\\
    &= \mathcal{F}( u_\mathcal{T}-w)- a_\mathcal{T}(u_\mathcal{T}, u_\mathcal{T}-w) - m_\mathcal{T}(u_\mathcal{T},u_\mathcal{T}-w).
\end{align*}  
\end{linenomath}
    Defining $e_\mathcal{T}:= u_\mathcal{T}-w$, we get
\begin{linenomath}
 \begin{align}\label{BoundStab1}
    \gamma S_\mathcal{T}(u_\mathcal{T}, e_\mathcal{T})&=\sum_{E\in \mathcal{T}}\left\{\int_E f e_\mathcal{T}- \int_E A_E  \Pi^0_{k-1}\left(\nabla u_\mathcal{T}\right)\;\Pi^0_{k-1}\left(\nabla e_\mathcal{T}\right) - \int_E c_E \Pi^0_{k} u_\mathcal{T}\; \Pi^0_{k} e_\mathcal{T}\right\}.
\end{align}
\end{linenomath}   
We notice that 
\begin{equation}\label{BoundStab2}
\begin{split}
    \int_E A_E \Pi^0_{k-1}\left (\nabla u_\mathcal{T}\right)\Pi^0_{k-1} \left(\nabla e_\mathcal{T}\right) &= \int_E \Pi^0_{k-1} \left(A_E \Pi^0_{k-1}\nabla u_\mathcal{T}\right)\nabla e_\mathcal{T} \\
    & =\int_E ( \Pi^0_{k-1} -I) (A_E \Pi^0_{k-1}\nabla u_\mathcal{T}) \nabla e_\mathcal{T} + \int_E A_E \Pi^0_{k-1}\nabla u_\mathcal{T} \nabla e_\mathcal{T}
\end{split}
\end{equation}
and 
\begin{linenomath}
\begin{align}\label{BoundStab3}
    \int_E  c_E \Pi^0_k u_\mathcal{T}\;\Pi^0_k e_\mathcal{T} = \int_E \Pi^0_k (c_E \Pi^0_k u_\mathcal{T}) \;e_\mathcal{T} = \int_E (\Pi^0_k - I)(c_E \Pi^0_k u_\mathcal{T}) e_\mathcal{T} + \int_E c_E (\Pi^0_k u_\mathcal{T}) e_\mathcal{T}.
\end{align}
\end{linenomath}
By substituting \eqref{BoundStab2} and \eqref{BoundStab3} into \eqref{BoundStab1}, it results
\begin{linenomath}
\begin{align*}
    \gamma S_\mathcal{T}(u_\mathcal{T},u_\mathcal{T}) &=
    \\&=\sum_{E\in \mathcal{T}}\int_E \left(f +  \nabla \cdot \left(A_E \Pi^0_{k-1} \nabla u_\mathcal{T}\right) - c_E \Pi^0_k u_\mathcal{T}\right) e_\mathcal{T}-\sum_{E\in \mathcal{T}}\int_{\partial E}\bm{n}\cdot \nabla \left(A_E \Pi^0_{k-1} \nabla u_\mathcal{T}\right)e_\mathcal{T}\\
    & \qquad +\sum_{E\in \mathcal{T}} \int_E (I-\Pi^0_{k-1})(A_E \Pi^0_{k-1}\nabla u_\mathcal{T})\nabla e_\mathcal{T} + \sum_{E\in \mathcal{T}} \int_E (I-\Pi^0_{k})(c_E \Pi^0_k u_\mathcal{T})\; e_\mathcal{T}\\
    &\le\sum_{E \in \mathcal{T}}h_E||r_\mathcal{T}(E; u_\mathcal{T}, \mathcal{D})||_{0,E} h^{-1}_E||e_\mathcal{T}||_{0,E} + \frac{1}{2} \sum_{e \in \mathcal{E}} h_e^{1/2}||j_\mathcal{T}(e; u_\mathcal{T}, \mathcal{D})||_{0,e}h^{-1/2}_e ||e_\mathcal{T}||_{0,e}\\
    & \qquad +\sum_{E\in \mathcal{T}} ||(I-\Pi^0_k)(A_E \Pi^0_{k} \nabla u_\mathcal{T})||_{0, E} ||\nabla e_\mathcal{T}||_{0, E} 
    + \sum_{E \in \mathcal{T}}h_E||(I-\Pi^0_k)c_E \Pi^0_k u_\mathcal{T}||_{0,E} h^{-1}_E ||e_\mathcal{T}||_{0, E}
    \\&\le\sum_{E \in \mathcal{T}}h_E||r_\mathcal{T}(E; u_\mathcal{T}, \mathcal{D})||_{0,E} h^{-1}_E||e_\mathcal{T}||_{0,E} + \frac{1}{2} \sum_{e \in \mathcal{E}} h_e^{1/2}||j_\mathcal{T}(e; u_\mathcal{T}, \mathcal{D})||_{0,e}h^{-1/2}_e ||e_\mathcal{T}||_{0,e}\\
    & \qquad +C_{\text{inv}}\sum_{E\in \mathcal{T}}\Psi_A(E; u_\mathcal{T}, \mathcal{D}) h^{-1}_E||e_\mathcal{T}||_{0, E} 
    + \sum_{E \in \mathcal{T}}\Psi_c(E; u_\mathcal{T}, \mathcal{D}) h^{-1}_E||e_\mathcal{T}||_{0, E}    .
\end{align*}
\end{linenomath}
With the same strategy used in \cite{AVEMstabfree}, for any $\delta>0$, we get
\begin{linenomath}
\begin{align*}
    \gamma S_\mathcal{T}(u_\mathcal{T},u_\mathcal{T})\le \frac{1}{2 \delta} \left( \eta^2_\mathcal{T}(u_\mathcal{T}, \mathcal{D}) + \Psi^2_\mathcal{T}(u_\mathcal{T}, \mathcal{D})\right)+ \frac{\delta}{2} \Phi_\mathcal{T}(e_\mathcal{T}),
\end{align*}
\end{linenomath}
where
\begin{linenomath}
\begin{align*}
    \Phi_\mathcal{T}(e_\mathcal{T})= \sum_{E \in \mathcal{T}}\Big \{ \max\{C^2_{\text{inv}},1\}h^{-2}_E|| e_\mesh ||_{0,E} + \frac{1}{2} \sum_{e \in \mathcal{E}} h^{-1}_E|| e_\mesh||_{0,e} \Big\}.
\end{align*}
\end{linenomath}
Posing now $w = \mathcal{I}^0_\mathcal{T}u_\mathcal{T}$ and applying Proposition \ref{ScaledPoincare}, we get
\begin{linenomath}
\begin{align*}
     \Phi_\mathcal{T}(u_\mathcal{T} -\mathcal{I}^0_\mathcal{T}u_\mathcal{T} ) \lesssim |u_\mathcal{T} -\mathcal{I}^0_\mathcal{T}u_\mathcal{T}|^2_{1,\Omega} \,,
\end{align*}
\end{linenomath}
whereas Proposition \ref{Prop:InterpolationOperators} yields
\begin{linenomath}
\begin{align*}
    |u_\mesh -\mathcal{I}^0_\mathcal{T} u_\mesh |^2_{1,\Omega}\lesssim  |u_\mesh -\mathcal{I}_\mathcal{T} u_\mesh|^2_{1,\Omega}\simeq S_\mathcal{T}(u_\mesh,u_\mesh),
\end{align*}
\end{linenomath}
so we obtain
\begin{linenomath}
\begin{align*}
    \gamma^2 S_\mathcal{T}(u_\mathcal{T},u_\mathcal{T})\le C_B  \left(\eta^2_\mathcal{T}(u_\mathcal{T}, \mathcal{D}) + \Psi^2_\mathcal{T}(u_\mathcal{T}, \mathcal{D})\right),
\end{align*}
\end{linenomath}
for a suitable constant $C_B>0$.
\end{proof}
Combining Propositions \ref{upperBound}  and \ref{bSt}, we arrive at the following key result. 
\begin{cor}[stabilization-free a posteriori error upper bound]\label{cor:stabFreeUpperBound} It holds
\begin{linenomath}
\begin{align*}
  |u - u_\mathcal{T}|^2_{1,\Omega}\le C_{U_1} \eta^2_\mathcal{T} (u_\mathcal{T}, \mathcal{D}) + C_{U_2}\Psi^2_\mathcal{T}(u_\mathcal{T}, \mathcal{D}),
\end{align*}
\end{linenomath}
where $C_{U_1} = C_\text{apost} \left(\frac{C_B}{\gamma^2} +1\right)$ and $C_{U_2} = C_\text{apost}\frac{C_B}{\gamma^2} $.
\end{cor}

\section{The effect of a mesh refinement}\label{sec:EffectOfAMeshRefinement}

In view of the convergence analysis of the adaptive algorithm {\tt GALERKIN}, in
this section we analyse the effect of refining the partition $\mesh$ by applying one or more newest-vertex bisections to some of its elements.  
Specifically, in Sect. \ref{sec:red-ref} we prove that the residual estimator \eqref{etaLocT} is reduced by a fixed fraction (up to an addend proportional to the stabilization term) when the element $E$ is split into two elements by one bisection. We prove a similar result for the inconsistency term estimator \eqref{def:PsiAc}, provided a suitable number of bisections is applied to $E$. Next, in Sect. \ref{sec:quasi-orthog} we establish a quasi-orthogonality property in the energy norm between the solutions on two partitions, one being a refinement of the other.

\subsection{Reduction of estimators under refinement}\label{sec:red-ref}

  Let us consider an element $E$ in $\mathcal{T}$ which is bisected into elements $E_1$ and $E_2$; the refined partition containing these two elements will be denoted by $\meshs$. Given $v \in \Vmesh$, we notice that $v$ is known on $\partial E$, and in particular at the new vertex of $E_1$ and $E_2$ produced by the bisection. Denoting by $e = E_1\cap E_2$ the new edge, we associate a function $v_*\in\mathbb{V}_{\mathcal{T}_*}$ to $v$ such that $v_*|_{\partial E} = v|_{\partial E}$, $v_*|_{e}\in \mathbb{P}_1(e)$, and $\bm{\mu}_p(E_i,v_*) =  \bm{\mu}_p(E, v)$ for all $0 \leq p \leq k-2$ and for $i=1,2$.
In the following we will write $v$ instead of $v_*$ when no confusion arises.

\subsubsection{\underline{The residual estimator}}

 Let $\eta_\mathcal{T}(E;v, \mathcal{D})$ be  defined in \eqref{etaLocT} and $\eta_{\mathcal{T}_*}(E;v, \mathcal{D})$ be the sum of the local residual estimators on the two newly formed elements , defined as follows:
\begin{linenomath}
\begin{align*} 
    \eta^2_{\mathcal{T}_*}(E; v, \mathcal{D}):= \sum^2_{i=1} \eta^2_{\mathcal{T}_*}(E_i; v, \mathcal{D})= \sum^2_{i=1} \left \{ h^2_{E_i} ||r_\mathcal{T}(E_i; v, \mathcal{D})||^2_{0,E_i} + \frac{1}{2} \sum_{e \in \mathcal{E}_{E_i}}h_{E_i} ||j_\mathcal{T}(e; v, \mathcal{D})||^2_{0,e} \right\},
\end{align*}
\end{linenomath}
where we recall that $h_ {E_i}= \frac{1}{\sqrt{2}}h_E$, $i=1,2$
We notice that, since $\mathcal{D}$ does not change under refinement, the functions $f_{E_i} = f_E|_{E_i}$, $c_{E_i} = c_E|_{E_i}$ and $A_{E_i} = A_E|_{E_i}$ will be denoted again by $f_E$, $c_E$ and $A_E$, respectively.

\begin{lemma}[local residual estimator reduction]\label{Lemma:LocalEstimatorReduction}
There exist constants $\mu_r \in (0,1)$ and $c_{er,1}>0$ such that for any $v \in \Vmesh$
\begin{linenomath}
\begin{align*} 
\eta_{\mathcal{T}_*}(E;v, \mathcal{D})\le \mu_r \;\eta_{\mathcal{T}}(E;v,\mathcal{D}) + c_{er,1} S^{1/2}_{\mathcal{T}(E)}(v,v),
\end{align*}
\end{linenomath}
where $S_{\mathcal{T}(E)}(v,v):= \sum_{E' \in \mathcal{T}(E)}S_{E'}(v,v)$ with $\mathcal{T}(E):= \{E'\in \mathcal{T}: \mathcal{E}_E \cap \mathcal{E}_{E'} \not=\emptyset\}$.
\end{lemma}
\begin{proof}
Recalling the definition \eqref{ResLoc}, we have the following residuals
\begin{linenomath}
\begin{align*}
    &r_E:= f_E + \nabla \cdot \left(A_E  \Pi_{k-1,E}^0 \nabla v\right) -c_E  \Pi_{k,E}^0 v \,,\\
    &r_{E_i}:= f_E + \nabla \cdot \left(A_E \Pi_{k-1,E_i}^0 \nabla v \right) - c_E \Pi_{k,E_i}^0 v \,.
\end{align*}
\end{linenomath}
Writing
$    r_{E_i} = r_E- \nabla \cdot\left( A_E  \Pi_{k-1,E}^0 \nabla v -  A_E \Pi_{k-1,E_i}^0 \nabla v\right)  + c_E \Pi_{k,E}^0 v - c_E \Pi_{k,E_i}^0 v $,
we get, for any $\epsilon>0$,
\begin{linenomath}
\begin{align*}
  \sum^2_{i=1} h^2_{E_i}||r_{E_i}||^2_{0, E_i} &\le    \sum^2_{i=1} h^2_{E_i}(1+ \epsilon)||r_E||^2_{0, E_i}\\
   &+  2\sum^2_{i=1} h^2_{E_i}\left( 1+ \frac{1}{\epsilon}\right)||\nabla \cdot \left(  A_E \left(\Pi_{k-1,E}^0 \nabla v - \Pi_{k-1,E_i}^0\nabla v\right)\right)||^2_{0,E_i}\\
   &+  2\sum^2_{i=1} h^2_{E_i}\left( 1+ \frac{1}{\epsilon}\right)||  c_E \left(\Pi_{k,E}^0 v  - \Pi_{k,E_i}^0 v \right)||^2_{0,E_i}.
\end{align*}
\end{linenomath}
The second term can be bounded by using the inverse inequality and the minimality of $\Pi_{k-1,E_i}^0$ as follows:
\begin{linenomath}
\begin{align*}
    \sum^2_{i=1} h^2_{E_i}||\nabla \cdot \left( A_E \left(\Pi_{k-1,E}^0 \nabla v - \Pi_{k-1,E_i}^0 \nabla v\right)\right)||^2_{0,E_i} &\lesssim   \sum^2_{i=1} ||   \Pi_{k-1,E}^0 \nabla v - \Pi_{k-1,E_i}^0\nabla v||^2_{0,E_i}\\
    & \le 2 ||\nabla v -\Pi_{k-1,E}^0 \nabla v||^2_{0, E}+ 2 \sum_{i=1}^2 ||\nabla v - \Pi_{k-1,E_i}^0 \nabla v||^2_{0,E_i}\\
    &  \le 4 |\nabla v -\Pi_{k-1,E}^0 \nabla v||^2_{0, E}
      \lesssim | v -\mathcal{I}_E v|^2_{1, E}\lesssim S_E(v,v) \,,
\end{align*}
\end{linenomath}
while, for the last term, using the Poincar\'e inequality we have
\begin{linenomath}
\begin{align*}
     \sum^2_{i=1} h^2_{E_i}||  c_E \left(\Pi_{k,E}^0 v - \Pi_{k,E_i}^0 v\right) ||^2_{0,E_i} &\lesssim h_E^2  \sum^2_{i=1} || \Pi_{k,E}^0 v - \Pi_{k,E_i}^0  v||^2_{0,E_i}\\
     &\le  h^2_E ||v - \Pi_{k,E}^0 v  ||^2_{0,E} \lesssim  h^2_E  | v - \Pi_{k,E}^0 v  |^2_{1,E}  \lesssim h^2_E S_E(v,v).
\end{align*}
\end{linenomath}
Finally, taking an appropriate value of $\epsilon$ and setting  $\mu:= \frac{1+ \epsilon}{2}\in (0, 1)$ (for instance, if  $\epsilon = \frac{1}{2}$,  $\mu= \frac{3}{4}$) we get
\begin{linenomath}
\begin{align*}
\sum^2_{i=1} h^2_{E_i}||r_{E_i}||^2_{0, E_i} &\le \mu \; h^2_{E}||r_E||^2_{0, E}+ C(1+ h^2_E )S_E(v,v) \,,
\end{align*}
\end{linenomath}
where $C>0$ is a constant.

For the jump condition, we will essentially use the proof given in \cite[Lemma 5.2]{AVEMstabfree}. In particular, we write $j_{\mathcal{T}_*}(e;v)= j_\mathcal{T}(e;v)+\left(j_{\mathcal{T}_*}(e;v)- j_\mathcal{T}(e,v)\right)$ and for any $\epsilon >0$
\begin{linenomath}
    \begin{align*}
        \sum^2_{j=1} \sum_{e\in \mathcal{E}_{E_i}}h_{E_i}\normed{j_{\mathcal{T}_*}(e;v)}^2\le(1+\epsilon)\, T_1 +\left(1+\frac1\epsilon\right) T_2 \,,
    \end{align*}
\end{linenomath}
with $T_1:=\sum^2_{i=1}\sum_{e\in\mathcal{E}_{E_i}}h_{E_i}\normed{j_\mathcal{T}(e;v)}^2$ and $T_2 :=\sum^2_{i=1}\sum_{e\in\mathcal{E}_{E_i}}h_{E_i}\normed{j_{\mathcal{T}_*}(e;v)-j_\mathcal{T}(e;v)}^2$. On the new edge we notice that $j_\mathcal{T}(e; v)=0$, then, 
\begin{linenomath}
    \begin{align*}
        T_1\le \frac{1}{\sqrt{2}}\sum_{e\in\mathcal{E}_{E}}h_{E}\normed{j_\mathcal{T}(e;v)}^2.
    \end{align*}
\end{linenomath}
We now define $\mathcal{T}_*(E_i):=\{E'\in \mathcal{T}_*: \mathcal{E}_{E_i}\cap \mathcal{E}_{E'}\neq \emptyset\}$;  for any edge $e\in \mathcal{E}_{E_i}$, we denote by $E_{i,e}\in \mathcal{T}_*(E_i)$ the element such that $e = \partial E_i \cap \partial E_{i,e}$. Then, 
\begin{linenomath}
    \begin{align*}
        \normed{j_{\mathcal{T}_*}(e;v) - j_\mathcal{T}(e;v)}&= \normed{\, [\![ A (\Pi^0_{\mathcal{T}_*} -\Pi^0_{\mathcal{T}})\nabla v]\!] \, }\\&\le \normed{A_E (\Pi^0_{k-1,E_i}-\Pi^0_{k-1,E})\nabla v} + \normed{A_{\hat{E}_{i,e}} (\Pi^0_{k-1,E_{i,e}}-\Pi^0_{k-1,\hat{E}_{i,e}})\nabla v},
    \end{align*}
\end{linenomath}
where $\hat{E}_{i,e}$ indicates the parent of $E_{i,e}$. Using the trace inequality we have
\begin{linenomath}
    \begin{align*}
        T_2&\lesssim \sum^2_{i=1}\sum_{E'\in \mathcal{T}_*(E_i)}||(\Pi^0_{k-1,E'} -\Pi^0_{k-1\hat{E'}} )\nabla v||_{0,E'}^2\\&\lesssim \sum^2_{i=1}\sum_{E'\in \mathcal{T}_*(E_i)}\left(||\nabla v -\Pi^0_{k-1,E'} \nabla v||_{0,E'}^2 +||\nabla v -\Pi^0_{k-1,\hat{E'}} \nabla v||_{0,E'}^2\right)
    \end{align*}
\end{linenomath}
Using now the minimality property of $\Pi^0_{k-1,E'}$ and $\Pi^0_{k-1,\hat{E}'}$, we easily get as above
\begin{linenomath}
    \begin{align*}
        T_2\le \sum_{E'\in \mathcal{T}(E)}||\nabla(v- \mathcal{I}_{E'}v)||^2_{0,E'} \lesssim \sum_{E'\in \mathcal{T}(E)} S_{E'}(v,v),
    \end{align*}
\end{linenomath}
which, for a sufficiently small $\epsilon$, concludes the proof.
\end{proof}

From this Lemma and the Lipschitz continuity of the residual estimator with respect to the argument $v$ (whose proof  is independent of the used polynomial degree, so we refer to \cite[ Lemma 5.3]{AVEMstabfree}), we immediately deduce the following result.

\begin{prop}[residual estimator reduction on refined elements]\label{prop: reduction}
There exist constants ${ \mu}_r \in (0,1)$, ${c_{er,1}}>0$ and ${ c_{er,2}}>0$ independent of $\mesh$ such that for any $v \in \Vmesh$ and $w \in \Vmeshs$, and any element $E \in \mesh$ which is split into two children $E_1,E_2 \in \meshs$, one has
\begin{linenomath}
\begin{equation}\label{eq:residual-comp00}
\etameshs(E;w,\data) \leq { \mu}_r \ \etamesh(E; v,\data) + { c_{er,1}} \, S^{1/2}_{\mesh(E)}(v,v) + { c_{er,2}} \, |v-w|_{1,\mesh(E)} \,.
\end{equation}
\end{linenomath}
\end{prop}

\subsubsection{\underline{The virtual inconsistency estimator}}\label{sec:virtual-reduction}

Given $v \in \Vmesh$ and $E \in \mesh$, consider the two virtual inconsistency terms $\Psi_{\mesh,A}(E,v,\data)$ and $\Psi_{\mesh,c}(E,v,\data)$ introduced in \eqref{eq:def-inconsistency}. When $E$ is bisected into $E_1$ and $E_2$, the term $\Psi_{\mesh,c}(E,v,\data)$ is reduced by a factor $\mu_c <1$ up to an addend proportional to the stabilization term, i.e., there exists $c_{vi,c} >0$ such that
\begin{linenomath}
\begin{equation}\label{eq:incon-psi-c}
\left( \sum_{i=1}^2 \Psi_{\meshs,c}(E_i,v,\data)^2\right)^{1/2}\!\!\!  \leq \mu_c \, \Psi_{\mesh,c}(E,v,\data) + c_{vi,c} \SE(v,v)^{1/2} \,.
\end{equation}
\end{linenomath}
This stems from the presence of the factor $h_E$ in front of the norm $||\left (I-\Pi^0_k)(c_E \Pi^0_k v \right)||_{0,E}$, with an argument similar to the one used in the proof of Lemma \ref{Lemma:LocalEstimatorReduction}.

Due to the lack of the factor $h_E$, a reduction result similar to \eqref{eq:incon-psi-c} does not hold for $\Psi_{\mesh,c}(E,v,\data)$. Indeed, since $A_E \Pi_{k-1,E}^0 \nabla v \in \mathbb{P}_{2k-2}(E)$, one may ask whether a constant $\mu<1$ esists such that
\begin{linenomath}
\begin{equation}\label{eq:refinement-numerics-1}
\sum_{i=1}^2 \Vert (I-\Pi_{k-1,E_i}^0 ) q \Vert_{0,E_i}^2 \le \mu^2  \Vert (I-\Pi_{k-1,E}^0 ) q \Vert_{0,E}^2 \qquad \forall q \in \mathbb{P}_{2k-2}(E) \,.
\end{equation}
\end{linenomath}
Unfortunately, the answer is no, as it can be seen numerically, working on the reference element $\hat{E}$ by affinity and identifying $\mu^2$ as the largest eigenvalue of a generalized eigenvalue problem. However, the same numerics indicates that if $\hat{E}$ is split into $2^m$ triangles of equal area by $m$ successive levels of uniform bisections, then $\mu^2$ becomes $<1$ for $m$ large enough, as seen in Table \ref{table:results}.
\begin{table}[t!]
    \centering
    \begin{tabular}{c|c|c}
             & $m=1$ &$m=2$\\
        \hline
        $k=2$ & 1.0000 &0.3153\\
        \hline
        $k=3$ & 1.0000 &0.6648
    \end{tabular}
    \caption{Value of $\mu^2$ in \eqref{eq:refinement-numerics-1} for different values of the polynomial degree $k$ and the level of refinement $m$}
    \label{table:results}
\end{table}
This is indeed predicted by the following result.
\begin{lemma}\label{lemma:reduction-numerics}
Let $E \in \mesh$. For any polynomial degree $k \geq 1$ there exists a minimal $m \in \mathbb{N}$ and a constant $\mu=\mu_m<1$ independent of $E$ such that, if $E$ is partitioned into $2^m$ elements $E_i$ of equal area by $m$ levels of uniform newest vertex bisection, it holds
\begin{linenomath}
\begin{equation}\label{eq:refinement-numerics-2}
\sum_{i=1}^{2^m} \Vert (I-\Pi_{k-1,E_i}^0 ) q \Vert_{0,E_i}^2 \le \mu^2  \Vert (I-\Pi_{k-1,E}^0 ) q \Vert_{0,E}^2 \qquad \forall q \in \mathbb{P}_{2k-2}(E) \,.
\end{equation}
\end{linenomath}
\end{lemma}
\begin{proof}
Since by construction $h_{E_i} = 2^{-m/2} h_E$, classical approximation results give
\begin{linenomath}
$$
\sum_{i=1}^{2^m} \Vert (I-\Pi_{k-1,E_i}^0 ) q \Vert_{0,E_i}^2 \le C_k 2^{-m} h_E^2 \vert q \vert_{1,E}^2
$$
\end{linenomath}
for some constant $C_k$ depending on $k$. Replacing $q$ by $q- \Pi_{k-1,E}^0  q$ leaves the left-hand side unchanged, whereas 
on the right-hand side an inverse inequality yields
\begin{linenomath}
$$
\sum_{i=1}^{2^m} \Vert (I-\Pi_{k-1,E_i}^0 ) q \Vert_{0,E_i}^2 \le C_k C_{\text{inv}, k}2^{-m}  \Vert q- \Pi_{k-1,E}^0  q \Vert_{0,E}^2 \,.
$$
\end{linenomath}
One concludes taking as $m$ the smallest integer such that $\mu_m^2 := C_k C_{\text{inv}, k}2^{-m} <1$.
\end{proof}

Based on these results, let $\meshs^{\! m}$ be a refinement of $\mesh$ in which the element $E$ has undergone $m$ levels of uniform refinements by newest vertex bisection, and has been replaced by $2^m$ subelements $E_i$. 
Given $v \in \Vmesh$, let us set
\begin{linenomath}
$$
  \Psi^2_{\meshs^{\!\! m},A}(E; v, \mathcal{D}) = \sum_{i=1}^{2^m}  \Vert(I-\Pi^0_{E_i, k-1}) (A_E\Pi^0_{E_i,k-1}\nabla v) \Vert_{0,E_i}^2 \,. 
$$
\end{linenomath}
\begin{lemma}\label{lemma:incon-psi-A}
There exist constants $\rho_A <1$ and $c_{vi,A} > 0$ such that for any $v \in \mathbb{V}_{E,k}$
\begin{linenomath}
$$
\Psi_{ \meshs^{\!\! m},A}(E; v, \mathcal{D}) \leq \rho_A  \Psi_{\mesh,A}(E; v, \mathcal{D}) + c_{vi,A} \SE^{1/2}(v,v) \,.
$$
\end{linenomath}
\end{lemma}
\begin{proof}
Write
\begin{linenomath}
\begin{equation*}
\begin{split}
\Vert(I-\Pi^0_{E_i, k-1}) (A_E\Pi^0_{E_i,k-1}\nabla v) \Vert_{0,E_i} &\leq \Vert(I-\Pi^0_{E_i, k-1}) (A_E\Pi^0_{E,k-1}\nabla v) \Vert_{0,E_i} \\
& \qquad +\Vert A_E (\Pi^0_{E_i,k-1}\nabla v-\Pi^0_{E,k-1}\nabla v) \Vert_{0,E_i} \,,
\end{split}
\end{equation*}
\end{linenomath}
sum over $i$, and conclude using \eqref{eq:refinement-numerics-2} and the usual arguments based on the minimality of the $L^2$-orthogonal projections.
\end{proof}

Let us set
\begin{linenomath}
$$
\Psi_{ \meshs^{\!\! m}}^2(E,v,\data) := \Psi_{  \meshs^{\!\! m},A}^2(E; v, \mathcal{D}) + \Psi_{  \meshs^{\!\! m},c}^2(E; v, \mathcal{D})
$$
\end{linenomath}
with
\begin{linenomath}
$$
\Psi_{  \meshs^{\!\! m},c}^2(E; v, \mathcal{D}) = \sum_{i=1}^{2^m} h_{E_i}^2 \Vert(I-\Pi^0_{E_i, k}) (c_E\Pi^0_{E_i,k}v) \Vert_{0,E_i}^2 \,. 
$$
\end{linenomath}
Applying a bound similar to \eqref{eq:incon-psi-c} to the successive level of refinements, we arrive at the following result.
\begin{lemma}\label{lemma:incon-psi}
There exist constants $\mu_{vi} <1$ and $c_{vi,1} > 0$ such that for any $v \in \mathbb{V}_{E,k}$
\begin{linenomath}
$$
\Psi_{  \meshs^{\!\! m}}(E; v, \mathcal{D}) \leq \mu_{vi} \,   \Psi_{\mesh}(E; v, \mathcal{D}) + c_{vi,1}\, \SE^{1/2}(v,v) \,.
$$
\end{linenomath}
\end{lemma}

Combining this estimate with the Lipschitz continuity property of the virtual inconsistency estimator, we obtain the following result.

\begin{prop}[virtual inconsistency estimator reduction on refined elements]\label{prop: reduction-vi}
There exist constants ${ \mu}_{vi} \in (0,1)$, ${c_{vi,1}}>0$ and ${ c_{vi,2}}>0$ independent of $\mesh$ such that for any $v \in \Vmesh$ and $w \in \V_{ \meshs^{\!\! m}} $, and any element $E \in \mesh$ which is split into $2^m$ children $E_i \in \V_{ \meshs^{\!\! m}} $, one has
\begin{linenomath}
\begin{equation}\label{eq:residual-comp00-vi}
\Psi_{  \meshs^{\!\! m}}(E; w, \mathcal{D})  \leq { \mu}_{vi} \, \Psi_{  \mesh}(E; v, \mathcal{D})   + { c_{vi,1}} \, S^{1/2}_{E}(v,v) + { c_{vi,2}} \, |v-w|_{1,E} \,.
\end{equation}
\end{linenomath}
\end{prop}

\subsection{Quasi-orthogonality property}\label{sec:quasi-orthog}

Let $u_{\mathcal{T}_*} \in \mathbb{V}_{\mathcal{T}_*}$ be  the solution of Problem \eqref{Discrete_Variazional_Problem} on the refined mesh $\mathcal{T}_*$. Hereafter we establish relations between the two energy errors $ \EnergyNorm{u - u_\mathcal{T}}$ and $\EnergyNorm{u - u_{\mathcal{T}_*}}$. 
The first result follows from Proposition \ref{Prop:InterpolationOperators} and Lemma \ref{Lemma:QuasiOrtog}; the proof is independent of the used polynomial degree, so we refer to \cite[Proposition 5.7]{AVEMstabfree}.

\begin{prop}[comparison of the energy error under refinement]\label{prop:comp-err} 
For any $\delta \in (0,1]$ there exists a constant $C_E>0$ independent of $\mathcal{T}$ and $\delta$ such that
\begin{linenomath}
\begin{align*}
    \EnergyNorm{u - u_{\mathcal{T}_*}}^2\le (1 + \delta)\EnergyNorm{u - u_{\mathcal{T}}}^2 - \EnergyNorm{u_{\mathcal{T}_*} - u_\mathcal{T}}^2 + C_E \left( 1+ \frac{1}{\delta}\right)\left(S_\mathcal{T}(u_\mathcal{T},u_\mathcal{T}) + S_{\mathcal{T}_*}(u_{\mathcal{T}_*},u_{\mathcal{T}_*})   \right).
\end{align*}
\end{linenomath}
\end{prop}

Next result extends Corollary 5.8 in \cite{AVEMstabfree}.
\begin{prop}[quasi-orthogonality of energy errors without stabilization]\label{prop:quasiOrtWithoutStab}
Given any $\delta \in \left(0, \frac{1}{4}\right]$, there exists $\gamma_{\delta}>0$ such that for any $\gamma>\gamma_\delta$, it holds 
\begin{linenomath}
\begin{align*}
    \EnergyNorm{u- u_{\mathcal{T}_*}}^2 \le (1 + 4 \delta)\EnergyNorm{u- u_{\mathcal{T}}}^2 - \EnergyNorm{u_{\mathcal{T}_*}- u_{\mathcal{T}}}^2 + 2\delta \left(\Psi^2_\mathcal{T}(u_\mathcal{T}, \mathcal{D}) + \Psi^2_{\mathcal{T}_*}(u_{\mathcal{T}_*}, \mathcal{D}) \right) \,.
\end{align*}
\end{linenomath}
    \begin{proof}
        Let $e:= \EnergyNorm{u -u_\mathcal{T}}$, $e_*:= \EnergyNorm{u -u_{\mathcal{T}_*}}$, $S:= S_\mathcal{T}(u_\mathcal{T},u_\mathcal{T})$, $S_*:= S_{\mathcal{T}_*}(u_{\mathcal{T}_*},u_{\mathcal{T}_*})$, $\eta:=\eta_\mathcal{T}(u_\mathcal{T}, \mathcal{D})$,  $\Psi:=\Psi_\mathcal{T}(u_\mathcal{T}, \mathcal{D})$, $\Psi_*:=\Psi_{\mathcal{T}_*}(u_{\mathcal{T}_*}, \mathcal{D})$ and $E:= \EnergyNorm{u_\mathcal{T} - u_{\mathcal{T}_*}}$.
        From Corollary \ref{cor:globalLoweBound} and \eqref{bound_energynorm}, we get 
        \begin{linenomath}
        \begin{align*}
            \eta^2 \le \frac{S}{c_{apost}} + \frac{e^2}{c_{apost}\;c_\mathcal{B}} \,,
        \end{align*}
        \end{linenomath}
    while, from Proposition \ref{bSt},
    \begin{linenomath}
    \begin{align*}
        S\le \frac{C_B}{\gamma^2}\left(\eta^2 + \Psi^2\right).
    \end{align*}
    \end{linenomath}
Combining them, we have
\begin{linenomath}
\begin{align*}
    \left(1- \frac{C_B}{\gamma^2\; c_{apost}} \right)S\le \frac{C_B}{\gamma^2}\left( \frac{e^2}{c_{apost}\;c_\mathcal{B}} + \Psi^2\right).
\end{align*}
\end{linenomath}
Doing the same on $\mesh_*$ and defining \begin{linenomath}   
 \begin{align*}
    \overline{C}:= \left(1- \frac{C_B}{\; c_{apost}} \right)^{-1} C_B \max\left\{1, \frac{1}{c_{apost}\;c_\mathcal{B}}\right\}
\le \left(1- \frac{C_B}{\gamma^2\; c_{apost}} \right)^{-1} C_B \max\left\{1, \frac{1}{c_{apost}\;c_\mathcal{B}}\right\}\end{align*}\end{linenomath} provided $\gamma^2\ge1$, we get
\begin{linenomath}
\begin{align*}
    &S\le \frac{\overline{C}}{\gamma^2} \left( e^2 + \Psi^2\right), &S_*\le \frac{\overline{C}}{\gamma^2} \left( e_*^2 + \Psi^2_*\right).
\end{align*}
\end{linenomath}
Employing Proposition \ref{prop:comp-err}, we obtain
\begin{linenomath}
\begin{align*}
    e^2_* \le (1 + \delta)e^2 -E^2 + C_E \left( 1 + \frac{1}{\delta} \right)\frac{\overline{C}}{\gamma^2}(e^2 + e_*^2 + \Psi^2 + \Psi^2_*). 
\end{align*}
\end{linenomath}
If we define $D:= C_E\left(  1 + \frac{1}{\delta} \right)\, \overline{C}$,
\begin{linenomath}
\begin{align*}
\left( 1 - \frac{D}{\gamma^2}\right)  e^2_* \le   \left(1 + \delta + \frac{D}{\gamma^2}\right)e^2 -E^2 + \frac{D}{\gamma^2}( \Psi^2 + \Psi^2_*).
\end{align*}
\end{linenomath}
By choosing $\gamma$ such that
\begin{linenomath}
\begin{equation}\label{eq:condition-gamma-delta-0}
\frac{1}{\gamma^2}\le \frac{\delta}{D} \,,
\end{equation}
\end{linenomath}
we get
\begin{linenomath}
\begin{align*}
    (1- \delta)e^2_* \le (1+2 \delta)e^2 - E^2 + \delta(\Psi^2 +\Psi^2_*),
\end{align*}
\end{linenomath}
which concludes the proof by observing that $\frac{1 + 2 \delta}{1-\delta}\le 1 + 4 \delta$ and $\frac{\delta}{1-\delta}\le 2\delta$, when $\delta \le \frac{1}{4}$.
\end{proof}
\end{prop}

\section{The module {\tt  GALERKIN}}\label{Sec:GALERKIN}

Let us consider a $\Lambda$-admissible input mesh $\mathcal{T}_0$, a set of approximated data $\mathcal{D}$ which consist of piecewise polynomials of degree $k-1$ on $\mathcal{T}_0$, and a tolerance $\epsilon>0$. The call 
\begin{linenomath}
\begin{align*}
    [\mathcal{T}, u_{\mathcal{T}}] = {\tt GALERKIN}(\mathcal{T}_0, \mathcal{D},\epsilon)
\end{align*}
\end{linenomath}
produces a $\Lambda$-admissible refined mesh $\mathcal{T}$ and the Galerkin approximation $u_\mathcal{T}\in \mathbb{V}_\mathcal{T}$, such as 
\begin{linenomath}
\begin{align*}
    \EnergyNorm{u- u_\mathcal{T}}\le C_G \epsilon,
\end{align*}
\end{linenomath}
where $u$ is the solution of Problem \eqref{Variational_Problem} and $C_G=\sqrt{c^\mathcal{B} \max\left\{C_{U_1},C_{U_2}\right\}}$, with $c^\mathcal{B}$ is defined in \eqref{bound_energynorm} and $C_{U_1},C_{U_2}$ in Corollary \ref{cor:stabFreeUpperBound}. We obtain it by iterating the sequence
 \begin{linenomath}
    \begin{align*}
    {\tt SOLVE} \rightarrow {\tt ESTIMATE} \rightarrow {\tt MARK} \rightarrow {\tt REFINE} \,.
\end{align*}
\end{linenomath}
At each step, a $\Lambda-$admissible mesh $\mathcal{T}_j$ and the associated solution $u_j$ of the discrete Problem \eqref{Discrete_Variazional_Problem} are produced. The process stops when the condition $\eta^2_{\mathcal{T}_j}(u_j, \mathcal{D}) + \Psi^2_{\mesh_j}(u_j, \data) \le \epsilon^2$ is reached. 

In particular, the modules are defined as follows:
\begin{itemize}
    \item $[u_\mathcal{T}] =  {\tt SOLVE}(\mathcal{T}, \mathcal{D})$ produces the solution of Problem \eqref{Discrete_Variazional_Problem} with data $\mathcal{D}$;
    \item $[\{\eta_\mathcal{T}(\cdot; u_\mathcal{T}, \mathcal{D})\}, \{\Psi_\mathcal{T}(\cdot; u_\mathcal{T}, \mathcal{D})\}]= {\tt ESTIMATE}(\mathcal{T}, u_\mathcal{T})$ computes the local estimators on $\mathcal{T}$;
    \item $[\mathcal{M}]= {\tt MARK}(\mathcal{T}, \{\eta_{\mathcal{T}}(\cdot; u_\mathcal{T}, \mathcal{D}) \}, \{\Psi_\mathcal{T}(\cdot; u_\mathcal{T}, \mathcal{D})\}] ,\theta)$ implements the D\"orfler criterion \cite{Dorfler} and finds  an almost minimal set $\mathcal{M}$  of elements in $\mathcal{T}$ such that
 \begin{linenomath}
    \begin{align}\label{RefinepropM}
        &\theta \left( \eta^2_\mathcal{T}(u_\mathcal{T}, \mathcal{D}) +  \Psi^2_\mathcal{T}(u_\mathcal{T}, \mathcal{D}) \right) \le \sum_{E\in \mathcal{M}} \left( \eta^2_\mathcal{T}(E; u_\mathcal{T}, \mathcal{D}) + \Psi^2_\mathcal{T}(E;u_\mathcal{T}, \mathcal{D}) \right),
    \end{align} 
    \end{linenomath}   
    for a given parameter $\theta \in (0,1)$;
    \item $[\mathcal{T}_*] = {\sf REFINE}(\mathcal{T},\mathcal{M}, \Lambda)$ returns a $\Lambda$-admissible refined mesh obtained from $\mathcal{T}$ by suitable newest-vertex bisections of the elements in $\mathcal{M}$, and possibly of other elements to fullfil the $\Lambda$-admissibility condition.
\end{itemize}

It is worth adding some details about the procedure {\sf REFINE}. Let $E \in \mathcal{M}$ be an element marked for refinement. For simplicity, hereafter let us set $\eta:=\eta_\mathcal{T}(E; u_\mathcal{T}, \mathcal{D})$ and $\Psi:= \Psi_\mathcal{T}(E; u_\mathcal{T}, \mathcal{D})$. The refinement of $E$ is performed as follows:
\begin{itemize}
\item if $\eta \geq \Psi$, then $E$ is bisected once;
\item if $\eta < \Psi$, then $E$ is bisected $m$-times, where $m$ has been introduced in Sect. \ref{sec:virtual-reduction} (see Lemma \ref{lemma:reduction-numerics}).
\end{itemize}
Denote by ${\cal P}(E)$ the partition of $E$ so obtained, and set $\eta_*^2 := \sum_{E' \in {\cal P}(E)} \eta_{{\cal P}(E)}^2(E'; u_\mesh, \data)$ and $\Psi_*^2 := \sum_{E' \in {\cal P}(E)} \Psi_{{\cal P}(E)}^2(E'; u_\mesh, \data)$. Then, recalling Lemma \ref{Lemma:LocalEstimatorReduction} and Lemma \ref{lemma:incon-psi}, one gets when $\eta \geq\Psi$
\begin{linenomath}
$$
\eta_* + \Psi_* \leq \frac{\mu_r+1}2 (\eta + \Psi) + c \, S^{1/2}_{\mathcal{T}(E)}(\umesh,\umesh)\,.
$$ \end{linenomath}
Indeed, $\Psi$ can be written as $\Psi = \lambda \eta$ for a certain $\lambda\in[0,1]$ and
\begin{linenomath}
\begin{align*}
\eta_*+\Psi_*&\leq \mu_r\eta +\lambda \eta+ c \, S^{1/2}_{\mathcal{T}(E)}(\umesh,\umesh) = \frac{\mu_r + \lambda}{1+\lambda}(1+\lambda)\eta+ c \, S^{1/2}_{\mathcal{T}(E)}(\umesh,\umesh) \\&=\frac{(\mu_r + \lambda)}{1+\lambda}(\eta +\Psi)+ c \, S^{1/2}_{\mathcal{T}(E)}(\umesh,\umesh).
\end{align*}\end{linenomath}
In the case $\eta< \Psi$,
\begin{linenomath}
$$
\eta_* + \Psi_* \leq \max(\mu_r^m, \mu_{vi}) (\eta + \Psi) + c \, S^{1/2}_{\mathcal{T}(E)}(\umesh,\umesh) \,.
$$\end{linenomath}
In all cases, it holds
\begin{equation}\label{eq:post-refine}
\eta_* + \Psi_* \leq \max\Big(\frac{\mu_r+1}2, \mu_{vi}\Big) (\eta + \Psi) + c \, S^{1/2}_{\mathcal{T}(E)}(\umesh,\umesh) \,,
\end{equation}
which shows that in each marked element the sum of the two estimators is reduced under refinement, up to the stabilization term. Note that for values $k=2$ or $3$ of the polynomial degree of practical use, two bisections ($m=2$) are enough when $\eta< \Psi$. 

This refinement may create non-admissible hanging nodes, i.e., hanging nodes with global index larger than $\Lambda$. To remove them and guaranteee $\Lambda$-admissibility of $\meshs$, further refinements should be applied. For the realization of this technical part, we refer to Sect. 11.1 in \cite{AVEMConvergenceOptimality}.


The following section proves the convergence of the {\tt GALERKIN} algorithm.

\section{Convergence property of {\tt GALERKIN}}\label{sec:galerkin}

\begin{prop}[global estimators reduction]\label{prop:estimatorReduction}
Let $u_\mathcal{T} \in \mathbb{V}_\mathcal{T}$ be the solution of the discrete variational Problem \eqref{Discrete_Variazional_Problem}. There exist constants $\rho \in (0,1)$ and  $C_{ger,1}, C_{ger,2} >0$ independent of $\mathcal{T}$ such that, if $\mathcal{T}_*$ is the refinement of $\mathcal{T}$ obtained by applying {\tt REFINE}, one has for any $w \in \mathbb{V}_{\mathcal{T}_*}$
\begin{equation} 
\begin{split}
        \eta_{\mathcal{T}_*}^2(w,\mathcal{D}) +  \Psi_{\mathcal{T}_*}^2(w,\mathcal{D})  &\le \rho \left( \etamesh^2(u_\mathcal{T},\mathcal{D})+ \Psi_\mathcal{T}^2(u_\mathcal{T},\mathcal{D}) \right) \\
        &  \qquad \quad + \ C_{ger,1} S_{\mathcal{T}}(u_\mathcal{T},u_\mathcal{T}) + C_{ger,2} \abs{u_\mathcal{T}-w}_{1, \Omega}^2\, \label{eq:estimatorReduction-1} \,.
\end{split}        
\end{equation}
\end{prop}
\begin{proof} One can reach the conclusion e.g. as in \cite[proof of Proposition 5.5]{AVEMstabfree}, using the bound \eqref{eq:post-refine} in each element $E$ marked for refinement.
\end{proof}

\begin{thm}[contraction property of {\sf GALERKIN}]
Let $\mathcal{M}\subset \mathcal{T}$ be the set of the marked elements relative to the solution $u_\mathcal{T} \in \mathbb{V}_\mathcal{T}$ of the discrete variational Problem \eqref{Discrete_Variazional_Problem}. If $\mathcal{T}_*$ is the refinement of $\mathcal{T}$ obtained by applying {\sf REFINE}, then for $\gamma$ sufficiently large there exist $\alpha \in (0,1)$ and $\beta >0$, $\zeta>0$ such that

\begin{linenomath}
\begin{equation*}
\EnergyNorm{ u -u_{\mathcal{T}_*}}^2 + \beta\, \eta^2_{\mathcal{T}_*}(u_{\mathcal{T}_*}, \mathcal{D}) + \zeta \, \Psi^2_{\mathcal{T}_*}(u_\mathcal{T}, \mathcal{D})   \le\alpha \left( \EnergyNorm{ u -u_{\mathcal{T}}}^2 + \beta \eta^2_{\mathcal{T}}(u_{\mathcal{T}}, \mathcal{D}) + \zeta \Psi^2_{\mathcal{T}}(u_\mathcal{T}, \mathcal{D}) \right).
\end{equation*}
\end{linenomath}
\end{thm}
\begin{proof} 
    To simplify notation, we set again $e= \EnergyNorm{u -u_\mathcal{T}}$, $e_*= \EnergyNorm{u -u_{\mathcal{T}_*}}$, $S= S_\mathcal{T}(u_\mathcal{T},u_\mathcal{T})$, $S_*= S_{\mathcal{T}_*}(u_{\mathcal{T}_*},u_{\mathcal{T}_*})$, $\eta=\eta_\mathcal{T}(u_\mathcal{T}, \mathcal{D})$, $\eta=\eta_{\mathcal{T}_*}(u_{\mathcal{T}_*}, \mathcal{D})$, $\Psi=\Psi_\mathcal{T}(u_\mathcal{T}, \mathcal{D})$, $\Psi_*=\Psi_{\mathcal{T}_*}(u_{\mathcal{T}_*}, \mathcal{D})$ and $E= \EnergyNorm{u_\mathcal{T} - u_{\mathcal{T}_*}}$. From Proposition \ref{prop:quasiOrtWithoutStab}, 
    \begin{linenomath}   
$$
        e^2_* \le ( 1+ 4 \delta)e^2 -E^2 +2\delta (\Psi +\Psi_*),
$$ \end{linenomath}
whereas using Proposition \ref{prop:estimatorReduction} and Proposition \ref{bSt}, we get
\begin{linenomath}
$$
        \eta^2_* + \Psi^2_* \le \rho ( \eta^2 + \Psi^2) + C_{ger,1}S + \frac{C_{ger,2}}{c_\mathcal{B}} E^2 \leq \left(\rho + \frac{C_{ger,1} C_B}{\gamma^2} \right) ( \eta^2 + \Psi^2)  +\frac{C_{ger,2}}{c_\mathcal{B}} E^2 \,.
$$\end{linenomath}
Combining them, we get
\begin{linenomath}
\begin{equation*}
\begin{split}
     e^2_* + \beta \eta^2_* + \left(\beta - 2\delta \right) \Psi^2_*  &\leq ( 1+ 4 \delta)e^2 +\left(  \frac{\beta C_{ger,2}}{c_\mathcal{B}} -1\right)E^2 \\ 
     & + \beta \left(\rho + \frac{C_{ger,1} C_B}{\gamma^2} \right) \eta^2 + \beta \left(\rho + \frac{C_{ger,1} C_B}{\gamma^2} + \frac{2\delta}{\beta} \right)  \Psi^2 \,, 
 \end{split}    
\end{equation*}    \end{linenomath} 
which suggests choosing $\beta$ such that
\begin{equation}\label{eq:cond-beta}
\frac{\beta C_{ger,2}}{c_\mathcal{B}} = 1 \,.
\end{equation}
Next, we write
\begin{linenomath}
\begin{equation*}
\begin{split}
     e^2_* + \beta \eta^2_* + \left(\beta - 2\delta \right) \Psi^2_*  &\leq ( 1- \delta)e^2 + 5\delta \, e^2  \\
 &   + \beta \left(\rho + \frac{C_{ger,1} C_B}{\gamma^2} \right) \eta^2 + \beta \left(\rho + \frac{C_{ger,1} C_B}{\gamma^2} + \frac{2\delta}{\beta} \right)  \Psi^2 \,, 
 \end{split}    
\end{equation*}\end{linenomath}     
and we invoke Corollary \ref{cor:stabFreeUpperBound} to write
\begin{linenomath}
$$
e^2 \leq {c^\cB}C_{apost} \Big(1+\frac{C_B}{\gamma^2}\Big)\eta^2 + {c^\cB}C_{apost}  \frac{C_B}{\gamma^2} \Psi^2  \,,
$$\end{linenomath}
which gives
\begin{linenomath}
\begin{equation*}
\begin{split}
     e^2_* + \beta \eta^2_* + \left(\beta - 2\delta \right) \Psi^2_*  &\leq ( 1- \delta)e^2 + \beta \left(\rho + \frac{C_{ger,1} C_B}{\gamma^2} + \frac{5\delta}{\beta} {c^\cB}C_{apost} \Big(1+\frac{C_B}{\gamma^2} \Big)  \right) \eta^2 \\ 
 &  \qquad \qquad + \beta \left(\rho + \frac{C_{ger,1} C_B}{\gamma^2} + \frac{2\delta}{\beta} + \frac{5\delta}{\beta} {c^\cB}C_{apost} \frac{C_B}{\gamma^2}  \right)  \Psi^2 \,.
 \end{split}    
\end{equation*}  \end{linenomath}
We now choose $\gamma$ and $\delta$ such that
\begin{linenomath}
$$
\rho + \frac{C_{ger,1} C_B}{\gamma^2} + \frac{5\delta}{\beta} {c^\cB}C_{apost} \Big(1+\frac{C_B}{\gamma^2} \Big)  \leq \frac{1+\rho}2
$$\end{linenomath}
which holds true if 
\begin{equation}\label{eq:cond-gamma-delta-1}
\frac{C_{ger,1} C_B}{\gamma^2} \leq \frac{1-\rho}4  \qquad \text{and} \qquad \frac{5\delta}{\beta} {c^\cB}C_{apost}(1+C_B) \leq \frac{1-\rho}4 
\end{equation}
(recall that we already assumed $\gamma^2\geq 1$). Similarly, we choose $\gamma$ and $\delta$ such that
\begin{linenomath}
$$
\beta \left(\rho + \frac{C_{ger,1} C_B}{\gamma^2} + \frac{2\delta}{\beta} + \frac{5\delta}{\beta} {c^\cB}C_{apost} \frac{C_B}{\gamma^2}  \right)  \leq (\beta - 2\delta) \frac{1+\rho}2 \,,
$$\end{linenomath}
which holds true if $\gamma$ satisfies the first condition in \eqref{eq:cond-gamma-delta-1}, whereas $\delta$ satisfies
\begin{equation}\label{eq:cond-gamma-delta-2}
\left( 2 + 5 {c^\cB}C_{apost} {C_B} + \frac{1+\rho}{\beta} \right)\delta \leq   \frac{1-\rho}4 \,.
\end{equation}
This proves the result, if we define $\zeta :=\beta - 2\delta$, with $\beta$ defined by \eqref{eq:cond-beta} and $\delta < \frac{\beta}2$, and 
\begin{equation}\label{eq:cond-alpha}
\alpha := \min \left(1-\delta, \frac{1+\rho}2 \right)\,.
\end{equation}
The conditions on $\gamma$ and $\delta$ which lead to the desired estimate are given in \eqref{eq:condition-gamma-delta-0}, \eqref{eq:cond-gamma-delta-1} and \eqref{eq:cond-gamma-delta-2}.
\end{proof}

\section{Conclusions}\label{sec:Conclusions}
In this paper, we presented an adaptive VEM of order $k\ge2$ on nonconforming triangular meshes. In the analysis, the space $\mathbb{V}^0_\mesh$ of continuous, piecewise polynomials functions of degree $k$ on the triangulation $\mesh$  plays a fundamental role. Indeed, it is contained in the global VEM space, $\mathbb{V}^0_\mesh  \subseteq \mathbb{V}_\mesh$, and guarantees a quasi-orthogonality property for any refinement $\meshs$ of $\mesh$, since $\V^0_\mesh \subseteq \V^0_\meshs$. By pivoting on this space, we proved an a posteriori error estimate which does not contain the stabilization term appearing in the VEM discrete formulation. Consequently, we established the convergence of the adaptive VEM algorithm, by a contraction argument. 

Extensions of our work include:
\begin{itemize}
    \item the complexity and optimality analysis of the two step algorithm AVEM mentioned in the Introduction to account for non-polynomial data;
    \item the study of a variant of the adaptive algorithm in which the polynomial degree $k$ may take large values, in the spirit of a $p$-version;
    \item the treatment of more general polygonal meshes.
\end{itemize}


\section*{Acknowledgments }
The authors performed this research in the framework of the Italian MIUR Award
“Dipartimenti di Eccellenza 2018-2022” granted to the Department of Mathematical Sciences, Politecnico di Torino (CUP: E11G18000350001).  CC was partially supported by the Italian MIUR
through the PRIN grant 201752HKH8; DF thanks the INdAM-GNCS project “Metodi numerici per lo studio di strutture geometriche parametriche complesse” (CUP: E53C22001930001). The authors are members of the Italian INdAM-GNCS research group.

\section*{Conflict of interest}
The authors  declare no conflict of interest.

\end{document}